\documentclass[12pt]{amsart}
\usepackage{amsmath,amssymb,amsfonts,amsthm,amsopn,dsfont}
\usepackage{graphics}

\def\endproof{\hfill\vbox{\hrule height0.6pt\hbox{%
   \vrule height1.3ex width0.6pt\hskip0.8ex
   \vrule width0.6pt}\hrule height0.6pt
  }\medskip}

\newcommand{\tfa}{time-frequency analysis}

\newcommand{\ft}{Fourier transform}
\newcommand{\stft}{short-time Fourier transform}

\newcommand{\tf}{time-frequency}

\newcommand{\modsp}{modulation space}
\newcommand{\psdo}{pseudodifferential operator}

\newtheorem{tm}{Theorem}[section]
\newtheorem{lemma}[tm]{Lemma}

\newtheorem{theorem}{Theorem}[section]
\newtheorem{corollary}[theorem]{Corollary}

\newtheorem{proposition}[theorem]{Proposition}
\newtheorem{remark}[theorem]{Remark}

\newcommand{\beqa}{\begin{eqnarray*}}
\newcommand{\eeqa}{\end{eqnarray*}}

\DeclareMathOperator*{\supp}{supp}

\newcommand{\field}[1]{\mathbb{#1}}
\newcommand{\bR}{\field{R}} 
\newcommand{\bZ}{\field{Z}} 
\def\la{\lambda}

 \def\cF{\mathcal{F}} 
 \def\cS{\mathcal{S}}

 \def\cG{\mathcal{G}}
 \def\cM{\mathcal{M}}

 \def\cC{\mathcal{C}}

\def\a{\aleph}

\def\vgf{V_gf}

\def\rd{\bR^d}

\def\rdd{{\bR^{2d}}}

\def\lrd{L^2(\rd)}
\def\lrdd{L^2(\rdd)}
\def\zd{\bZ^d}

\def\intrd{\int_{\rd}}

\def\R{\right)}

\def\<{\left<}
\def\>{\right>}

\def\inv{^{-1}}

\def\mv1{M_v^1}

\def\mpq{M^{p,q}}

\def\phas{(x,\o )}
\def\mn{(m,n)}
\def\mn'{(m',n')}

\def\o{\eta}
\def\a{\alpha}
\def\b{\beta}

\def\N{\mathbb{N}}
\def\R{\mathbb{R}}
\def\Ren{\mathbb{R}^d}
\def\Renn{\mathbb{R}^{2d}}

\def\sch{\mathcal{S}}

\def\Fur{\mathcal{F}}

\def\f{\varphi}

\def\Sn2{S_{2}(L^{2}(\Ren))}
\def\S1{S_{1}(L^{2}(\Ren))}
\def\sig00{\sigma_{0,0}}

\def\la{\langle}
\def\ra{\rangle}

\newcommand{\sg}{\mathbf{SG}}
\newcommand{\tjk}{A_{j,k}}
\newcommand{\ttjk}{\tilde{A}_{j,k}}
\begin{document}

\begin{abstract}
We consider a class of Fourier integral
operators, globally defined on
$\R^{d}$, with symbols and phases
satisfying product type estimates (the
so-called $\sg$ or scattering classes).
We prove a sharp continuity result for
such operators when acting on the
modulation spaces $M^p$. The minimal
loss of derivatives is shown to be
$d|1/2-1/p|$. This global perspective produces
 a loss of decay as well, given
by the same order. Strictly related,
striking examples of unboundedness on
$L^p$ spaces are presented.
\end{abstract}

\title[On the global boundedness of Fourier integral operators]{On the global boundedness of Fourier integral operators}\author{Elena Cordero, Fabio Nicola and Luigi Rodino}
\address{Department of Mathematics,
University of Torino, via Carlo Alberto
10, 10123 Torino, Italy}
\address{Dipartimento di Matematica,
Politecnico di Torino, corso Duca degli
Abruzzi 24, 10129 Torino, Italy}
\address{Department of Mathematics,
University of Torino, via Carlo Alberto
10, 10123 Torino, Italy}
\email{elena.cordero@unito.it}
\email{fabio.nicola@polito.it}
\email{luigi.rodino@unito.it}
\thanks{}
\maketitle

\section{Introduction}
The Fourier integral operators (FIOs)
of H\"ormander
(\cite{hormander0,hormander,treves}),
in a simplified local version, are
operators of the form:
\begin{equation}\label{FIO}
Af(x)=A_{\Phi,\sigma}f(x)=\int e^{2\pi
i\Phi(x,\eta)}
\sigma(x,\eta)\hat{f}(\eta)\,d\eta.
\end{equation}
Here the \ft\, of $f\in\cS(\rd)$ is
normalized to be ${\hat
  {f}}(\o)=\int
f(t)e^{-2\pi i t\o}dt$. The phase
function $\Phi\phas$ in \eqref{FIO} is
 assumed real-valued, smooth for $\eta\not=0$ and positively homogeneous
 of degree 1 with respect to $\o$; moreover,
  $\sigma\phas$ belongs to H\"ormander's symbol
  class $S^m_{1,0}$
   of order $m\in\R$:
\begin{equation}\label{symb0}
|\partial^\alpha_\eta\partial^\beta_x
\sigma(x,\eta)|\leq C_{\alpha,\beta}
\langle\eta\rangle^{m-|\alpha|}, \quad
\forall (x,\eta)\in\R^{2d},
\end{equation}
where
$\langle\eta\rangle=(1+|\o|^2)^{1/2}$.
The definition being local, or
localized in a compact manifold, that
is, the support of $\sigma\phas$ is
assumed to have compact projection on
the space of the $x$-variables, say
$\sigma(x,\eta)=0$ for $|x|\geq R$, for
a suitable $R>0$. Moreover,
$\sigma(x,\eta)$ is usually cut to zero
near $\o=0$, that is $\sigma(x,\eta)=0$
in the strip
\begin{equation}\label{strip}
\{\phas\in\rdd,\quad|\o|\leq 1\}.
                \end{equation}
This eliminates the discontinuity at
$\o=0$ of the phase function
$\Phi\phas$ without no practical effect
on the local behaviour of the operator
$A$, since the eliminated part
corresponds to a (locally) regularizing
operator.

Boundedness in $\lrd$ and $L^p(\rd)$ of
$A$ have been widely studied, see e.g.
\cite{phongstein97, sogge93,stein93}
and references therein. As basic
results, we know that under the
non-degeneracy condition
\begin{equation}\label{nondeg}
\left|{\rm det}\,
\left(\frac{\partial^2\Phi}{\partial
x_i\partial \eta_l}\Big|_{
(x,\eta)}\right)\right|>\delta>0,\quad
\forall (x,\eta)\in \R^{2d},
\end{equation}
the operator $A$ is $L^2$-bounded for
$m=0$, see \cite{hormander0}, as well
as $L^p$-bounded,
 $1<p<\infty$, if the order $m$ of $\sigma\phas$ is negative, satisfying
\begin{equation}\label{vecchio}
m\leq-(d-1)\left|\frac{1}{2}-\frac{1}{p}\right|,
\end{equation}
see \cite{seegersoggestein} and
references quoted there.\par
 The result
cannot be improved in general, as clear
from the Fourier integral operator
solving the Cauchy problem for the wave
equation in space-dimension $d$. See
\cite{ruz1,ruz2} for a precise
discussion of the sharpness of
\eqref{vecchio}, depending on the
singular support of the kernel of $A$.
According to \eqref{vecchio}, in the
one-dimensional case, the assumption
$m=0$ is sufficient to get
$L^p$-boundedness for any $p$,
$1<p<\infty$.\par In
\cite{cordero-nicola-rodino} we studied
the action of an operator $A$ as above
on the spaces $\Fur L^p$ of temperate
distributions whose Fourier transform
is in $L^p$ (with the norm $\|f\|_{\Fur
L^p}=\|\hat{f}\|_{L^p}$). There it was
shown that $A$ is bounded as an
operator $(\Fur L^p)_{comp}\to (\Fur
L^p)_{loc}$, $1\leq p\leq\infty$, if $
m\leq-d\left|\frac{1}{2}-\frac{1}{p}\right|$.
This is similar to \eqref{vecchio}, but
with the difference of one unit in the
dimension. Surprisingly, this threshold
was shown to be sharp in any dimension
$d\geq1$, even for phases linear with
respect to $\eta$; see
\cite{cordero-nicola-rodino} (or
Section 6 below) for the construction
of explicit counterexamples.\par In the
present paper we want to study the
global boundedness of Fourier integral
operators as in \eqref{FIO}. Namely, we
consider the case when the support of
$\sigma\phas$ is not compact with
respect to the space variable $x$. In
this direction, general
$L^2$-boundedness results can be found
in \cite{ruzhsugimoto}; to this paper
we address for references on previous
$L^2$-global results and for
motivations, mainly concerning
hyperbolic problems where
global-in-space information is needed.

As a preliminary step of our study,
 we call attention on the following striking,
 but seemingly unknown, example.
 In dimension $d=1$, consider
 \begin{equation}\label{FIOes}
Af(x)=\int_\R e^{2\pi i \f(\eta)x}
\sigma(x,\eta)\hat{f}(\eta)\,d\eta,
\end{equation}
where $\sigma\in S^{0}_{1,0}$, and
$\f:\R\to\R$ is a diffeomorphism,
 with $\f(\o)=\o$ for $|\o|\geq 1$ and whose
 restriction to $(-1,1)$ is non-linear.
 This can be regarded as a
 pseudodifferential operator
 with symbol $e^{2\pi i
 x(\f(\eta)-\eta)}\sigma(x,\eta)$,
 which satisfies the
 estimates
 in \eqref{symb0}, with $m=0$, for $x$ in
 bounded subsets of $\R^d$.
 Hence it is bounded as an
 operator $L^p\to L^p_{loc}$,
 $1<p<\infty$ (\cite[page 250]{stein93}).
 Naively, one may think that the uniform bounds
 \eqref{symb0} for $\sigma$, with $m=0$, grant global $L^p$-boundedness as well. Instead we have:

\begin{theorem}\label{mo} Let
$2<p<\infty$. Assume $\sigma\phas=1$ in
\eqref{FIOes}; then $A$ is not bounded
as an operator from $L^p(\R)$ to
$L^p(\R)$. More precisely, fix
$\sigma\phas=\la x\ra^{\tilde{m}},$
$\tilde{m}\in \R$, in \eqref{FIOes};
then, $A: L^p(\R)\to L^p(\R)$ is
bounded if and only if
\begin{equation}\label{sf}
\tilde{m}\leq-\left(\frac12-\frac1p\right).
\end{equation}
\end{theorem}
Observe that, if $\sigma\phas=1$,
\textit{microlocally} for $|\o|\geq 1$
the operator $A$ is the identity
operator.
 Hence the behaviour of $\Phi\phas$ in
 the strip \eqref{strip} is now crucial and we
  could as well take
\begin{equation}\label{newsy}
\sigma\phas=\la
x\ra^{\tilde{m}}G(\o),\quad
\mbox{with}\,\,\,G\in\cC_0^\infty(\R),\,\,\,
G(\o)=1\,\,\mbox{for}\,\,|\o|\leq 1,
\end{equation}
as symbol in \eqref{FIOes}, without
changing the conclusions.\\ In the
subsequent Proposition \ref{casolp} we
present similar examples in every
dimension $d\geq1$ and for every $1\leq
p\leq\infty$, obtaining the threshold
\begin{equation}\label{vecchio1}
\tilde{m}\leq-d\left|\frac{1}{2}-\frac{1}{p}\right|,
\end{equation}
which is the same as that for local
$\Fur L^p$ spaces; see also Coriasco
and Ruzhansky \cite{coriascoruz} for
other examples in this connection.\par

Results of global $L^p$-boundedness,
taking simultaneously account of
\eqref{vecchio} and \eqref{vecchio1},
are given in the forthcoming paper
\cite{coriascoruz}.

The approach here will be different.
Namely, inspired by our previous papers
\cite{fio1,cordero-nicola-rodino}, we
replace $L^p$ by other function spaces,
the so-called modulation spaces $M^p$,
introduced by Feichtinger in \cite{F1},
which will allow us to restore a
symmetry between the thresholds
\eqref{vecchio} and \eqref{vecchio1}.
To be definite, let us first be precise
about the class of FIOs we consider,
and then recall the definition of
$M^p$.\par\medskip {\bf Global Fourier
integral operators.} We will be
concerned here with a class of FIOs
\eqref{FIO} with phase $\Phi$ and
symbol $\sigma$ chosen in the so-called
$\sg$ classes. Namely, keeping locally
the H\"{o}rmander's estimates
\eqref{symb0}, we shall introduce a
precise scale for the decay as
$x\to\infty$. The symbol $\sigma\in
\mathcal{C}^\infty(\R^{2d})$ is assumed
to belong to the class
 $\sg^{m_1,m_2}$ (the so-called class of global
 symbols, or scattering symbols, of order
 $(m_1,m_2)$), i.e.
\begin{equation}\label{symb}
|\partial^\alpha_\eta\partial^\beta_x
\sigma(x,\eta)|\leq C_{\alpha,\beta}
\langle\eta\rangle^{m_1-|\alpha|}\langle
x\rangle^{m_2-|\beta|}, \quad \forall
(x,\eta)\in\R^{2d},
\end{equation}
see, e.g., Cordes \cite{co2}, Parenti
\cite{parenti72}, Melrose
\cite{mel1,mel2}, Schrohe \cite{sc1},
Schulze \cite{sch1}. Note that the
classes $\sg^{m_1,m_2}$ are stable
under conjugation by Fourier transform,
namely $\cF^{-1} \sg^{m_1,m_2}\cF=
\sg^{m_2,m_1}$. Corresponding FIOs were
considered by Coriasco
\cite{coriasco,coriasco1,coriasco2},
Cappiello \cite{cappiello}, Cordes
\cite{cordes}, see also
Ruzhansky-Sugimoto \cite{ruzhsugimoto}
and references therein. The phase
function $\Phi(x,\eta)$ is real-valued
  and in the class $\sg^{1,1}$.
We also assume the non-degeneracy
condition \eqref{nondeg}. The operator
$A$ in \eqref{FIOes} is of this type,
having symbol $\sigma\in
\sg^{0,\tilde{m}}$ or, if $\sigma$ is
as in \eqref{newsy}, $\sigma\in
\sg^{-\infty,\tilde{m}}$. Locally, the
corresponding FIOs are of the type
\eqref{symb0}, with a somewhat more
general phase function; in particular,
for local $L^p$-boundedness the
threshold \eqref{vecchio} still holds
true. Global $L^2$-boundedness follows
from \cite{coriasco}; see also
\cite{ruzhsugimoto} for a more general
class of
 FIOs. Finally, we would like to address to
 the recent monography of Cordes
 \cite{cordes} for the role of $\sg$
  pseudodifferential operators and FIOs
  in Dirac's
  theory.\par\medskip
\textit{\bf Modulation spaces}. We
briefly recall the definition of the
modulation spaces $M^p$, $1\leq p\leq
\infty$, which are widely used in
time-frequency analysis (see
\cite{F1,book} and Section 2 for
definition and properties). In short,
we say that a temperate distribution
$f$ belongs to $M^p(\R^d)$ if its
short-time Fourier transform $V_g
f\phas$, defined in \eqref{STFT} below,
is in $L^p(\rdd)$, namely if
\begin{equation}\label{modulas}
\|f\|_{M^p}:=\|\|f(\cdot)\overline{g(\cdot-x)}\|_{\Fur
L^p}\|_{L^p_x}<\infty.
\end{equation}
 Here $g$ is a non-zero
(so-called window) function in
$\cS(\rd)$, which in \eqref{modulas} is
first translated and then multiplied by
$f$ to localize $f$ near any point $x$.
Changing $g\in \cS(\rd)$ produces
equivalent norms. The space
$\tilde{M}^{\infty}(\rd)$ is
 the closure of $\cS(\rd)$ in the
 $M^\infty$-norm. For heuristic purposes,
distributions in $M^p$ may be regarded
as functions which are locally in
${\Fur L^p}$ and decay at infinity like
functions in $L^p$ (see Lemma
\ref{lloc} below for a precise
statement). Among their properties, we
highlight their stability under \ft:
$\cF (M^p)=M^p$, $1\leq p\leq\infty$
(and $\cF
(\tilde{M}^\infty)=\tilde{M}^\infty$).

We may now state our result.
\begin{theorem}\label{maintheorem}
Let $\sigma\in \sg^{m_1,m_2}$ and
$\Phi\in\sg^{1,1}$ satisfying
\eqref{nondeg}. If
\begin{equation}\label{soglia}
m_1\leq-d\left|\frac{1}{2}-\frac{1}{p}\right|,\quad
m_2\leq-d\left|\frac{1}{2}-\frac{1}{p}\right|,
\end{equation}
then the corresponding FIO A, initially
defined on $\mathcal{S}(\R^d)$, extends
to a bounded operator on $M^{p}$,
whenever $1\leq p<\infty$. For
$p=\infty$, $A$ extends to a bounded
operator on $\tilde{M}^\infty$.
\end{theorem}
Both the bounds in \eqref{soglia} are
sharp. Namely, for any
$m_1>-d\displaystyle{\left|{1/2}-{1/p}\right|}$,
or
$m_2>-d\displaystyle{\left|{1/2}-{1/p}\right|}$,
there exists $A$ as in \eqref{FIO} with
$\sigma\in \sg^{m_1,-\infty}$,
$\sigma\in \sg^{-\infty,m_2}$,
respectively, ($\sigma$ being compactly
supported with respect to $x$ and
$\eta$ respectively) which is not
bounded on $M^p$.\par Let us compare
Theorem \ref{maintheorem} with our
preceeding results
\cite{fio1,cordero-nicola-rodino}. In
\cite{fio1} we considered different
Fourier integral operators,
corresponding to operator solutions to
Schr\"odinger equations, basic example
of phase functions being quadratic
forms in the $x,\o$ variables. Such
operators were proved to be bounded on
$M^p$ without loss of derivatives,
i.e., for symbols $\sigma\phas$ of
order zero, see also
\cite{benyi,concetti-toft,cordero2}. In
\cite{cordero-nicola-rodino} we
considered local H\"{o}rmander's FIOs
and proved that they are $M^p$ bounded
with the \textit{sharp} loss of
regularity $-d|1/2-1/p|$, i.e. the same
of that of operators acting on local
$\Fur L^p$ spaces. This agrees with the
loss for $m_1$ in Theorem
\ref{maintheorem}. Moreover, in Theorem
\ref{maintheorem} a further
 loss of decay (that
for $m_2$) appears, which agrees with
that of the example in Theorem \ref{mo}
for the action on global $L^p$ spaces.
This circle of relationships is well
understood by means of the heuristic
interpretation, given above, of the
modulation spaces.
 Also, we
underline that the invariance under
Fourier conjugation of the modulation
spaces $M^p$ reveals them to be an
appropriate functional framework for
global FIOs (this is the insight the
reader will catch from the proofs in
the sequel). \par Finally we observe
that these results should extend to the
more general class of global FIOs
considered in \cite{ruzhsugimoto}; we
plan to devote a subsequent paper to
this investigation.

\medskip
The paper is organized as follows. In
Section \ref{section2} the definitions
and basic properties of the modulation
spaces $M^p$ are recalled. Section
\ref{section3} contains a review of
$\sg$ FIOs and a boundedness result for
a class of $\sg$ FIOs whose phases have
bounded second derivatives (Proposition
\ref{pro2}). In Section \ref{sectionb}
we prove boundedness results on
modulation spaces for $\sg$
pseudodifferential operators. Section
\ref{section5} is devoted to the proof
of Theorem \ref{maintheorem}. Finally,
Section \ref{sharp} exhibits the
optimality of Theorem \ref{maintheorem}
and shows the negative results for
operators acting on $L^p$ spaces,
extending the example \eqref{FIOes} in
Theorem \ref{mo} above.

\vskip0.5truecm

\par

\textbf{Notation.} We define
$|x|^2=x\cdot x$, for $x\in\Ren$, where
$x\cdot y=xy$ is the scalar product on
$\Ren$. The space of smooth functions
with compact support is denoted by
$\cC_0^\infty(\rd)$, the Schwartz class
is $\sch(\Ren)$, the space of tempered
distributions $\sch'(\Ren)$.
 Translation and modulation operators ({\it time and frequency shifts}) are defined, respectively, by
$$ T_xf(t)=f(t-x)\quad{\rm and}\quad M_{\o}f(t)= e^{2\pi i \o
 t}f(t).$$
We have the formulas $(T_xf)\hat{} =
M_{-x}{\hat {f}}$, $(M_{\o}f)\hat{}
=T_{\o}{\hat {f}}$, and
$M_{\o}T_x=e^{2\pi i x\o}T_xM_{\o}$.
The inner product of two functions
$f,g\in \lrd$ is $\la f, g \ra=\intrd
f(t) \overline{g(t)}\,dt$, and its
extension to $\cS'\times\cS$ will be
also denoted by $\la \cdot, \cdot \ra$.

Given a weight function $\mu$ defined
on some lattice $\Lambda$, the spaces
$\ell^{p,q}_\mu$ are the Banach spaces
of sequences $\{a_{m,n}\}_{m,n}$,
$(m,n)\in\Lambda$, such that
$$\|a_{m,n}\|_{\ell^{p,q}_\mu}:=\left(\sum_{n}\left(\sum_{m}
|a_{m,n}|^p\mu(m,n)^p\right)^{q/p}\right)^{1/q}<\infty
$$
(with obvious changes when $p=\infty$
or $q=\infty)$.

The notation $A\lesssim B$ means $A\leq
c B$ for a suitable constant $c>0$,
whereas $A \asymp B$ means $c\inv A
\leq B \leq c A$, for some $c\geq 1$.
The symbol $B_1 \hookrightarrow B_2$
denotes the continuous embedding of the
space $B_1$ into $B_2$.

\section{Preliminary results on Time-Frequency methods}\label{section2}

First we summarize some concepts and
tools of \tfa, now available in
textbooks~\cite{folland89,book}. We
also recall some results from
\cite{fio1,cordero-nicola-rodino}.

\subsection{Modulation spaces}
The \stft\ (STFT) of a distribution
 $f\in\sch'(\Ren)$ with
respect to a non-zero window
$g\in\sch(\Ren)$ is
\begin{equation}\label{STFT} V_gf(x,\o)=\la
f,M_\o T_x g\ra =\int_{\Ren}
 f(t)\, {\overline {g(t-x)}} \, e^{-2\pi i\o t}\,dt.
 \end{equation}

The STFT $\vgf $ is defined on many
pairs of Banach spaces. For instance,
it maps $L^2(\rd ) \times L^2(\rd )$
into $\lrdd $ and
$\sch(\Ren)\times\sch(\Ren)$ into
$\sch(\Renn)$. Furthermore, it can be
extended to a map from
$\sch'(\Ren)\times\sch(\Ren)$ into
$\sch'(\Renn)$.

Recall the inversion formula for the
STFT (see e.g. (\cite[Corollary
3.2.3]{book}): if $\|g\|_{L^2}=1$ and,
for example, $u\in{L^2}(\R^d)$, it
turns out
\begin{equation}\label{treduetre}
u=\int_{\R^{2d}} V_g u(y,\o)M_\o T_y
g\, dy\,d\o.
\end{equation}


\label{modspdef}

The modulation space norms are a
measure of the joint time-frequency
distribution of $f\in \sch '$. For
their basic properties we refer, for
instance, to \cite[Ch.~11-13]{book} and
the original literature quoted there.

For the quantitative description of
decay and regularity properties, we use
 weight
functions on the \tf\ plane. In the
sequel $v$ will always be a continuous,
positive, even, submultiplicative
weight function (in short, a
submultiplicative weight), i.e.,
$v(0)=1$, $v(z) = v(-z)$, and $
v(z_1+z_2)\leq v(z_1)v(z_2)$, for all
$z, z_1,z_2\in\Renn.$ Associated to
every submultiplicative weight we
consider the class of so-called {\it
  v-moderate} weights $\cM _v$. A positive, even
weight function $\mu \neq 0$ everywhere
on $\Renn$ belongs to $\cM _v$ if it
satisfies the condition
$$
 \mu(z_1+z_2)\leq Cv(z_1)\mu(z_2) \quad \forall z_1,z_2\in\Renn \, .
$$
 We note that this definition implies that
$\frac{1}{v} \lesssim \mu \lesssim v $
and that $1/\mu \in \cM _v$.

 By abuse of notation, we denote product weights
  $v_{s_1,s_2}(x,\eta)=\langle
  x
\rangle ^{s_2}\langle\eta \rangle
^{s_1}$, $s_1,s_2\in\R$ (the indices'
order follows that of the
$\sg^{m_1,m_2}$-classes). Note that
$v_{s_1,s_2}$ is submultiplicative only
if $s_1,s_2\geq 0$.

Given a non-zero window
$g\in\sch(\Ren)$, a moderate weight
$m\in \cM _v$ and $1\leq p,q\leq
\infty$, the {\it
  modulation space} $M^{p,q}_{\mu}(\Ren)$ consists of all tempered
distributions $f\in\sch'(\Ren)$ such
that $V_gf\in L^{p,q}_{\mu}(\Renn )$
(weighted mixed-norm spaces). The norm
on $M^{p,q}_{\mu}$ is
$$
\|f\|_{M^{p,q}_{\mu}}=\|V_gf\|_{L^{p,q}_{\mu}}=\left(\int_{\Ren}
  \left(\int_{\Ren}|V_gf(x,\o)|^p{\mu}(x,\o)^p\,
    dx\right)^{q/p}d\o\right)^{1/p} \,
    ,
$$
with obvious changes if $p=\infty$ or
$q=\infty$. If $p=q$, we write
$M^p_{\mu}$ instead of $M^{p,p}_{\mu}$,
and if ${\mu}(z)\equiv 1$ on $\Renn$,
then we write $M^{p,q}$ and $M^p$ for
$M^{p,q}_{\mu}$ and $M^{p,p}_{\mu}$.

Then $\mpq_{\mu}(\Ren )$ is a Banach
space whose definition is independent
of the choice of the window $g$.
 Moreover,
if ${\mu}\in\cM_v$ and $g \in M^1_{v}
\setminus \{0\}$, then $\|V_gf
\|_{L^{p,q}_{\mu}}$ is an equivalent
norm for $M^{p,q}_{\mu}(\Ren)$ (see
\cite[Thm.~11.3.7]{book}). Roughly
speaking, a weight in $\o $ regulates
the smoothness of $f \in \mpq_{\mu}$,
whereas a weight in $x$ regulates the
decay at infinity.\par Denote by
$\tilde{{M}}^{p,q}_{\mu}$ the closure
of the Schwartz class in
$M^{p,q}_{\mu}$. We have
$\tilde{{M}}^{p,q}_{\mu}=M^{p,q}_{\mu}$
if $p<\infty$ and $q<\infty$ and the
duality property for \modsp s can be
stated as follows: if $1\leq p,q \leq
\infty $ and $p',q'$ are
 the conjugate exponents, then
 $(\tilde{{M}}^{p,q}_{\mu})^*=\tilde{{M}}^{p',q'}_{1/{\mu}}$.

For simplicity, we shall write
$M^{p,q}_{s_1,s_2}$ for
$M^{p,q}_{v_{s_1,s_2}}$ and, similarly,
for the spaces
$\tilde{{M}}^{p,q}_{s_1,s_2}$.

\par
 We also recall from \cite{Feichtinger-grochenig89}
  the following useful interpolation
  relations:\par
 {\it If $0<\theta<1$, $1\leq
p_1,p_2\leq \infty$, $1\leq
p\leq\infty$, and
 $s,\tilde{s},s_1,\tilde{s}_1,s_2,\tilde{s}_2\in\R$
 satisfy
 $1/p=(1-\theta)/p_1+\theta/p_2$,
  $s=(1-\theta)s_1+\theta s_2$,
  $\tilde{s}=(1-\theta)\tilde{s}_1
  +\theta \tilde{s}_2$, then}
\begin{equation}\label{interpmp}
(\tilde{M}^{p_1}_{s_1,\tilde{s}_1},
\tilde{M}^{p_2}_{s_2,\tilde{s}_2})_{[\theta]}=
\tilde{M}^p_{s,\tilde{s}}.
\end{equation}

For tempered distributions compactly
supported either in time or in
frequency, the $M^{p,q}$-norm is
equivalent to
 the $\cF L^q$-norm or $L^p$-norm, respectively.
  This result is well-known
  (\cite{fe89-1,feichtinger90}).
    See also \cite{kasso07}
   and \cite{cordero-sharp} for a proof.
\begin{lemma}\label{lloc} Let $1\leq p,q\leq
\infty.$\\
$(i)$ For every $u\in
\mathcal{S}'(\rd)$, supported in a
compact set $K\subset \rd$, we have
$u\in M^{p,q}\Leftrightarrow u\in \Fur
L^q$, and
\begin{equation}\label{loc}
C_K^{-1} \|u\|_{M^{p,q}}\leq \|u\|_{\cF
L^q}\leq C_K \|u\|_{M^{p,q}},
\end{equation}
where $C_K>0$ depends only on
 $K$.\\
 $(ii)$ For every $u\in \mathcal{S}'(\rd)$,
whose Fourier transform is supported in
a compact set $K\subset \rd$, we have
$u\in M^{p,q}\Leftrightarrow u\in L^p$,
and
\begin{equation}\label{loc2}
C_K^{-1} \|u\|_{M^{p,q}}\leq \|u\|_{
L^p}\leq C_K \|u\|_{M^{p,q}},
\end{equation}
where $C_K>0$ depends only on
 $K$.\\
\end{lemma}

In order to state the dilation
properties for modulation spaces, we
introduce the indices:
$$ \mu_1(p)=\begin{cases}-1/{p^\prime} & \quad {\mbox{if}} \quad 1\leq p\leq 2,\\
 -1/p & \quad {\mbox{if}} \quad 2\leq p\leq
 \infty,
 \end{cases}
 $$
and
$$ \mu_2(p)=\begin{cases}-1/p & \quad {\mbox{if}} \quad 1\leq p\leq 2,\\
 -1/{p^\prime} & \quad {\mbox{if}} \quad 2\leq p\leq \infty.\\
 \end{cases}
 $$
For $\lambda>0$, we define the dilation
operator $U_\lambda f(x)=f(\lambda x)$.
Then, the dilation properties of $M^p$
are as follows (see \cite[Theorem
3.1]{sugimototomita}).

\begin{theorem}\label{dilprop}
We have:
 $(i)$ For $\lambda\geq 1 $,
$$\| U_\lambda f\|_{M^p}\lesssim \lambda^{d\mu_1(p)}
\|f\|_{M^p},\quad\forall\, f\in
M^p(\rd).
$$
$(ii)$ For $0<\lambda\leq 1 $,
$$\| U_\lambda f\|_{M^p}\lesssim \lambda^{d\mu_2(p)}
\|f\|_{M^p},\quad\forall \, f\in
M^p(\rd).
$$
\end{theorem}
\noindent These dilation estimates are
sharp, as discussed in
\cite{sugimototomita}, see also
\cite{cordero2}.\par We also need the
following result.
\begin{lemma}\label{le41}
Let $\chi$ be a smooth function
supported where $B_0^{-1}\leq|\eta|\leq
B_0$, for some $B_0>0$.\par $(a)$ For
every $u\in\cS(\rd)$,
\[
\sum_{j=1}^\infty\|\chi(2^{-j}
D)u\|_{M^1}\lesssim \|u\|_{M^1},
\]
where $\chi(2^{-j} D)u=
\cF^{-1}[\chi(2^{-j}\cdot)\hat{u}]$.\par
$(b)$ For every $u\in\cS(\rd)$,
\[
\sum_{j=1}^\infty\|\chi(2^{-j}
\cdot)u\|_{M^1}\lesssim \|u\|_{M^1}.
\]
\end{lemma}
\begin{proof}
Part $(a)$ was proved in \cite[Lemma
5.1]{cordero-nicola-rodino}, whereas
part $(b)$ follows from $(a)$, since
the Fourier transform defines an
automorphisms of any $M^p$.
\end{proof}\par
Finally we recall the following result.
 \begin{lemma}\label{lemma2}
$(a)$ For $k\geq0$, let
$f_k\in\cS(\rd)$ satisfy ${\rm
supp}\,\hat{f}_0\subset B_2(0)$ and
\[
{\rm
\supp}\,\hat{f}_k\subset\{\eta\in\R^d:\
2^{k-1}\leq|\eta|\leq 2^{k+1}\},\quad
k\geq1.
\]
Then, if the sequence $f_k$ is bounded
in $M^\infty§(\rd)$, the series
$\sum_{k=0}^\infty {f}_k$ converges in
${M}^\infty(\R^d)$ and
\begin{equation}\label{b0}
\|\sum_{k=0}^\infty
f_k\|_{M^\infty}\lesssim\sup_{k\geq0}\|f_k\|_{M^\infty}.
\end{equation}
\indent $(b)$ For $k\geq0$, let
$f_k\in\cS(\rd)$ satisfy ${\rm
supp}\,{f}_0\subset B_2(0)$ and
\[
{\rm \supp}\,{f}_k\subset\{x\in\R^d:\
2^{k-1}\leq|x|\leq 2^{k+1}\},\quad
k\geq1.
\]
Then, if the sequence $f_k$ is bounded
in $M^\infty(\rd)$, the series
$\sum_{k=0}^\infty {f}_k$ converges in
${M}^\infty(\R^d)$ and \eqref{b0} holds
true.
 \end{lemma}
 \begin{proof}
Part $(a)$ was proved in \cite[Lemma
5.2]{cordero-nicola-rodino}, whereas
part $(b)$ follows from $(a)$, again
since the Fourier transform defines an
automorphisms of any $M^p$.
 \end{proof}

 \vspace{3 mm}
\subsection{Gabor frames}
Fix a function $g\in\lrd$ and a lattice
$\Lambda =\a \zd\times \b\zd$, for
$\a,\b>0$. For $(k,n)\in \Lambda$,
define $g_{k,n}:=M_{n}T_{k}g$. The set
of time-frequency shifts
$\cG(g,\a,\b)=\{g_{k,n}, (k,n)\in
\Lambda\}$ is called Gabor system.
Associated to $\cG(g,\a,\b)$ we define
the coefficient operator $C_g$, which
maps functions to sequences as follows:
\begin{equation}\label{analop}
    (C_gf)_{k,n}=(C_g^{\a,\b}f)_{k,n}:=\la
    f,g_{k,n}\ra,\quad (k,n)\in \Lambda,
\end{equation}
the synthesis operator
\[
D_g{c}=D_g^{\a,\b} c =\sum_{(k,n)\in
\Lambda} c_{k,n} M_{n}T_{k}g,\quad
c=\{c_{k,n}\}_{(k,n)\in \Lambda}
\]
and the Gabor frame operator
\begin{equation}\label{Gaborop}
    S_g f=S_g^{\a,\b}f:=D_g S_g f=\sum_{(k,n)\in
    \Lambda}\la f,g_{k,n}\ra g_{k,n}.
\end{equation}

The set $\cG(g,\a,\b)$ is called a
Gabor frame for the Hilbert space
$\lrd$ if $S_g$ is a bounded and
 invertible operator on $\lrd$. Equivalently, $C_g$ is bounded from $\lrd$ to $\ell^2(\a\zd\times\b\zd)$ with closed
  range, i.e., $\|f\|_{L^2}\asymp\|C_g f\|_{\ell^2}.$ If $\cG(g,\a,\b)$ is a Gabor frame for $\lrd$,
  then the so-called \emph{dual window}
   $\gamma=S_g^{-1} g$ is well-defined and the
   set $\cG(\gamma,\a,\b)$ is a frame (the so-called
   canonical dual frame of $\cG(g,\a,\b)$). Every $f\in \lrd$ possesses the frame expansion
\begin{equation}\label{frame}
    f= \sum_{(k,n)\in \Lambda}\la f,g_{k,n}\ra
     \gamma_{k,n}= \sum_{(k,n)\in \Lambda}\la f,
     \gamma_{k,n}\ra g_{k,n}
\end{equation}
with unconditional convergence in
$\lrd$, and norm equivalence:
$$\|f\|_{L^2}\asymp \|C_g f\|_{\ell^2}\asymp
\|C_\gamma f\|_{\ell^2}.
$$
This result is contained in
\cite[Proposition 5.2.1]{book}. In
particular, if $\gamma=g$ and
$\|g\|_{L^2}=1$ the frame is called
\emph{normalized tight} Gabor frame and
the expansion \eqref{frame} reduces to
\begin{equation}\label{parsevalframe}
    f= \sum_{(k,n)\in \Lambda}\la f,g_{k,n}\ra g_{k,n}.
\end{equation}
If we ask for more regularity on the
window $g$, then the previous result
can be extended to suitable Banach
spaces, as shown below
\cite{fg97jfa,GL01}.
\begin{theorem}\label{teomod}
Let $\mu\in\cM_v$, $\cG(g,\a,\b)$ be a
normalized tight Gabor frame for
$\lrd$, with lattice
$\Lambda=\a\zd\times\b \zd$,
 and $g\in M^1_v$. Define
 $\tilde{\mu}=\mu_{|_{\Lambda}}$.\\
(i) For every $1\leq p,q\leq\infty$,
$C_g: M^{p,q}_\mu\to
\ell^{p,q}_{\tilde{\mu}}$ and $D_g:
\ell^{p,q}_{\tilde{\mu}}\to
M^{p,q}_\mu$ countinuously and, if
$f\in\mpq_\mu,$ then the Gabor
expansions \eqref{parsevalframe}
converge unconditionally in $\mpq_\mu$
for $1\leq p,q<\infty$
and all weight $\mu$, and weak$^\ast$-$M^\infty_\mu$ unconditionally if $p=\infty$ or $q=\infty$.\\
(ii) The following norms are equivalent on $\mpq_\mu$:\\
\begin{equation}\label{framexp}
\|f\|_{\mpq_\mu}\asymp \|C_g f
\|_{\ell^{p,q}_{\tilde {\mu}}}.
\end{equation}
\end{theorem}
We also establish the following
property (\cite[Theorem 2.3]{fio1}).
Denote by
$\tilde{\ell}^{p,q}_{\tilde{\mu}}$ the
closure of the space of eventually zero
sequences in
$\ell^{p,q}_{\tilde{\mu}}$. Hence
$\tilde{\ell}^{p,q}_{\tilde{\mu}}={\ell}^{p,q}_{\tilde{\mu}}$
if $p<\infty$ and $q<\infty$.
\begin{theorem}\label{teomod2}
Under the assumptions of Theorem
\ref{teomod}, for every $1\leq
p,q\leq\infty$ the operator $C_g$ is
continuous from $\tilde{{M}}^{p,q}_\mu$
into
$\tilde{\ell}^{p,q}_{\tilde{\mu}}$,
whereas the operator $D_g$ is
continuous from
$\tilde{\ell}^{p,q}_{\tilde{\mu}}$ into
$\tilde{{M}}^{p,q}_\mu$.
\end{theorem}

\section{Preliminary results on
FIOs}\label{section3} The general
theory of SG FIOs was developed by
Coriasco in \cite{coriasco}, see also
\cite{cordes,coriasco1,coriasco2}. In
this section we recall the main
properties needed in the sequel. We
also present a boundedness result, in
the spirit of
\cite{fio1,cordero-nicola-rodino}, for
FIOs with phases having bounded
derivatives of order $\geq2$.\par
\subsection{SG
Fourier integral operators} First of
all we observe that the calculus for SG
FIOs was developed in \cite{coriasco}
for phases $\Phi\in \sg^{1,1}$
satisfying the growth condition
\begin{equation}\label{crescita}
\langle\nabla_x
\Phi(x,\eta)\rangle\gtrsim\langle\eta\rangle,\quad
\langle\nabla_\eta
\Phi(x,\eta)\rangle\gtrsim\langle
x\rangle.
\end{equation}
These conditions are a consequence of
our hypotheses,
 namely $\Phi\in \sg^{1,1}$ and
 \eqref{nondeg}. More precisely, it
 follows from the estimates
 from the mixed second
 derivatives, namely
\begin{equation}\label{crescita2}
|\partial^\alpha_\eta\partial^\beta_x
\Phi(x,\eta)|\leq
C_{\alpha,\beta},\quad
|\alpha|=|\beta|=1,\
\forall(x,\eta)\in\R^{2d},
\end{equation}
combined with \eqref{nondeg} and the
formula for
 the Jacobian of the inverse function, that
 Hadamard's
 global inverse function
 theorem (\cite{krantz}[Theorem
 6.2.4]{}) applies to the maps
 $x\longmapsto \nabla_\eta
\Phi(x,\eta)$ and
$\eta\longmapsto\nabla_x\Phi(x,\eta)$,
which are therefore globally
invertible. Moreover these maps have
bounded Jacobians uniformly with
respect to $\eta$ and $x$ respectively,
so that they are globally Lipschitz
continuous, uniformly with respect to
$\eta$ and $x$ respectively. The same
holds for their inverses, which proves
\eqref{crescita}.\par The first
important results is the following
formula for the composition of a $\sg$
pseudodifferential operator, namely an
operator of the form
\begin{equation}\label{pseudoformula}
p(x,D)u= \int e^{2\pi ix\eta}
p(x,\eta)\hat{f}(\eta)\,d\eta,
\end{equation}
with a symbol $p\in \sg^{t_1,t_2}$,
 and a
 $\sg$
FIO $A=A_{\Phi,\sigma}$ as in
\eqref{FIO} (\cite[Theorem
7]{coriasco}; for the case of
H\"ormander's symbol classes see
\cite{hormander,mascarello-rodino,treves}).

 First we recall that a regularizing operator
is a pseudodifferential operator
$R=r(x,D)$ with symbol $r(x,\eta)$ in
the Schwartz space
$\mathcal{S}(\R^{2d})$ (equivalently,
an operator with kernel in
$\mathcal{S}(\R^{2d})$, which maps
$\mathcal{S}'(\R^{d})$ into $
\mathcal{S}(\R^{d})$).
\begin{theorem}\label{composition}
Let the symbol $\sigma$ and the phase
$\Phi$ satisfy the assumptions in the
Introduction. Let $p(x,\eta)$ be a
symbol in $\sg^{t_1,t_2}$. Then,
\[
p(x,D)A=S+R,
\]
where $S$ is a FIO with the same phase
$\Phi$ and symbols $s(x,\eta)$ in the
class $\sg^{m_1+t_1,m_2+t_2}$,
satisfying
\[
{\rm supp}\,s\subset {\rm
supp}\,\sigma\cap\{(x,\eta)\in\R^{2d}:\
(x,\nabla_x\Phi(x,\eta))\in{\rm
supp}\,p\},
\]
and $R$ is a regularizing operator.\par
Moreover, the symbol estimates
satisfied by $s$ and the seminorm
estimates of $r$ in the Schwartz space
are uniform when $\sigma$ and $p$ vary
in bounded subsets of $\sg^{m_1,m_2}$
and $\sg^{t_1,t_2}$ respectively.\par
Also, the symbol $s$ has the following
asymptotic expansion:
\begin{equation}\label{asintotica}
s(x,\eta)\sim\sum_{\alpha\in\mathbb{Z}^d_+}
\frac{1}{\alpha!}
\partial^\alpha_\eta
p(x,\nabla_x\Phi(x,\eta))D_y^\alpha
[e^{i\psi(x,y,\eta)}
\sigma(y,\eta)]_{y=x},
\end{equation}
where
\[
\psi(x,y,\eta)=\Phi(y,\eta)-\Phi(x,\eta)-\langle
y-x,\nabla_x\Phi(x,\eta)\rangle,
\]
and, as usual,
$D^\alpha_y=(-i)^{|\alpha|}\partial^\alpha_y$.
\end{theorem}
The meaning of the above asymptotic
expansion is that the difference
between $s(x,\eta)$ and the partial sum
over $|\alpha|<N$ is a symbol in
$\sg^{m_1+t_1-N,m_2+t_2-N}$. However,
in the next sections only the first
part of the statement will be used.\par
Similarly, one also has the following
formula for the composition in the
reverse order (\cite[Theorem
8]{coriasco}). To state it, we
introduce the notation
$^t\!b(x,\eta):=b(\eta,x)$ for a
function $b(x,\eta)$ in $\R^{2d}$.
\begin{theorem}\label{composition2}
Let the symbol $\sigma$ and the phase
$\Phi$ satisfy the assumptions in the
Introduction. Let $p(x,\eta)$ be a
symbol in $\sg^{t_1,t_2}$. Then,
\[
Ap(x,D)=S+R,
\]
where $S$ is a FIO with the same phase
$\Phi$ and symbols $s(x,\eta)$ in the
class $\sg^{m_1+t_1,m_2+t_2}$,
satisfying
\[
{\rm supp}\,s\subset {\rm
supp}\,\sigma\cap\{(x,\eta)\in\R^{2d}:\
(\nabla_\eta\Phi(x,\eta),\eta)\in{\rm
supp}\,p\},
\]
and $R$ is a regularizing operator.\par
Moreover, the symbol estimates
satisfied by $s$ and the seminorm
estimates of $r$ in the Schwartz space
are uniform when $\sigma$ and $p$ vary
in bounded subsets of $\sg^{m_1,m_2}$
and $\sg^{t_1,t_2}$ respectively.\par
Also, the transpose symbol
$^t\!s(x,\eta)$ admits the asymptotic
expansion in \eqref{asintotica}, with
$p,\sigma$ and $\Phi$ replaced by
$^t\!p$, $^t\!\sigma$ and $^t\!\Phi$
respectively.
\end{theorem}
This latter result can be proved
combining Theorem \ref{composition}
with the following nice formula for the
transpose of $A=A_{\Phi,\sigma}$ with
respect to the pairing which extends
the integral $(u,v)\longmapsto\int uv$
(\cite[Proposition 9]{coriasco}):
\begin{equation}\label{congF}
^t\! A_{\Phi,\sigma}=\Fur \circ
A_{^t\!\Phi,^t\!\sigma}\circ\Fur^{-1}.
\end{equation}
Similarly, it is easily verified that
the $L^2$ formal adjoint of the FIO
$A_{\Phi,\sigma}$ is the operator
defined by
\begin{equation}\label{FIOB}
\widehat{Bf}(\o)=\widehat{B_{\Phi,\sigma}f}(\o)=\int
e^{-2\pi i\Phi(x,\eta)}
\overline{\sigma(x,\eta)}f(x)\,d x.
\end{equation}
namely,
\begin{equation}\label{aggiunto}
(A_{\Phi,\sigma})^\ast=B_{\Phi,\sigma}.
\end{equation}
In the sequel the operators of the type
\eqref{FIOB} will be called ``type II
FIOs". In contrast, operators of the
type \eqref{FIO} will be called ``type
I FIOs", or simply ``FIOs".\par
 Another important
result that will be used in the sequel
is the following one (\cite[Theorem
16]{coriasco}).
\begin{theorem}\label{elle2}
Let $\sigma\in\sg^{0,0}$ and
 $\Phi$ satisfying
the assumptions in the Introduction.
Then the corresponding FIO $A$,
initially defined on $\cS(\R^d)$,
extends to a bounded operator on
$L^2(\R^d)$.
\end{theorem}
The proof relies on the fact (cf.
\cite{hormander0}) that the composition
$A^\ast A$ is a pseudodifferential
operator with symbol in $\sg^{0,0}$;
therefore it is continuous on
$L^2(\R^d)$ (e.g., by \cite[Theorem
18.1.11]{hormander}). So $A$ is. This
result was generalized in
\cite{ruzhsugimoto} to FIOs with phases
satisfying weaker symbol estimates.

\subsection {FIOs with phases
having bounded derivatives of order
$\geq2$} We present here a boundedness
result of a class of FIOs whose phases
have bounded second derivatives (cf.
\cite[Theorem 4.1]{fio1} and
\cite[Proposition
3.3]{cordero-nicola-rodino}).
\begin{proposition}\label{pro2}
Consider a symbol $\sigma\in
C^\infty(\R^{2d})$ satisfying the
estimates
\begin{equation}\label{0i0tris}
|\partial^\alpha_\eta\partial^\beta_x\sigma(x,\eta)|
\leq C_{\alpha,\beta},\ \forall
(x,\eta)\in\R^{2d},
\end{equation}
 and a phase $\Phi\in
C^\infty(\R^{2d})$, satisfying
\begin{equation}\label{i0tris}
|\partial^\alpha_\eta\partial^\beta_x\Phi(x,\eta)|\leq
C_{\alpha,\beta}\quad {\rm for}\
|\alpha|+|\beta|\geq2,
\end{equation}
for $(x,\eta)$ in an
$\epsilon$-neighbourhood of the support
of $\sigma$, as well as
\begin{equation}\label{crescita2bis}
|\partial^\alpha_\eta\partial^\beta_x
\Phi(x,\eta)|\leq
C_{\alpha,\beta},\quad
|\alpha|=|\beta|=1,\
\forall(x,\eta)\in\R^{2d},
\end{equation}
 and
\begin{equation}\label{i1tris}
\left|{\rm det}\,
\left(\frac{\partial^2\Phi}{\partial
x_i\partial \eta_l}\Big|_{
(x,\eta)}\right)\right|>\delta>0,\quad
\forall (x,\eta)\in\R^{2d}.
\end{equation}
Then, for every $1\leq p\leq\infty$ it
turns out
\[
\|A u\|_{M^p}\leq C\|u\|_{M^p},\quad
\forall u\in\mathcal{S}(\R^d),
 \]
where the constant $C$ depends only on
$\delta, \epsilon$ and upper bounds for
a finite number
 of the
constants in \eqref{0i0tris},
\eqref{i0tris} and
\eqref{crescita2bis}.
\end{proposition}
\begin{proof} This is a
variant of \cite[Theorems 3.1,
4.1]{fio1} and \cite[Propositions 3.2,
3.3]{cordero-nicola-rodino}. For the
sake of completeness we outline the
proof. \par Let
$g,\gamma\in\mathcal{S}(\R^d)$,
$\|g\|_{L^2}= \|\gamma\|_{L^2}=1$, with
${\rm supp}\, \gamma\subset
B_{\epsilon/4}(0)$, ${\rm supp}\,
\hat{g}\subset B_{\epsilon/4}(0)$. Let
$u\in\mathcal{S}(\R^d)$. The inversion
formula \eqref{treduetre} for the STFT
gives
\[
V_\gamma
(Au)(y',\omega')=\int_{\R^{2d}}\langle
A(M_\omega T_y g), M_{\omega'}
T_{y'}\gamma\rangle V_g
u(y,\omega)dy\,d\omega.
\]
Hence, it suffices to prove that the
map $K_A$ defined by
\[
K_A G(y',\omega')=\int_{\R^{2d}}\langle
A(M_\omega T_y g), M_{\omega'}
T_{y'}\gamma\rangle
G(y,\omega)dy\,d\omega
\]
is continuous on $L^p(\R^{2d})$. By
Schur's test (see e.g. \cite[Lemma
6.2.1]{book}) we are reduced to proving
that its integral kernel
\[
K_A(y',\omega';y,\omega)=\langle
A(M_\omega T_y g),
M_{\omega'}T_{y'}\gamma\rangle
\]
satisfies
\begin{equation}\label{schur1}
K_A\in
L^\infty_{y',\omega'}(L^1_{y,\omega}),
\end{equation}
and
\begin{equation}\label{schur2}
K_A\in
L^\infty_{y,\omega}(L^1_{y',\omega'}).
\end{equation}

Now, in view of the hypothesis
\eqref{0i0tris} and \eqref{i0tris} we
can apply \cite[Proposition
3.2]{cordero-nicola-rodino}, that tells
us that for every $N\geq0$, there
exists a constant $C>0$ such that
\[
|\langle A (M_\omega T_y g),
M_{\omega'} T_{y'}\gamma\rangle|\leq C
\langle
\nabla_x\Phi(y',\omega)-\omega'\rangle^{-N}\langle
\nabla_\eta\Phi(y',\omega)-y\rangle^{-N}.
\]
The constant $C$ only depends on $N$,
$g,\gamma$, and on a finite number of
constants in \eqref{0i0tris} and
\eqref{i0tris}. \\
For $N>d$, $\intrd \langle
\nabla_\eta\Phi(y',\omega)-y\rangle^{-N}=\int
\langle y\rangle^{-N} \,dy<\infty$,
hence \eqref{schur1} will be proved if
we verify that there exists a constant
$C'>0$ such that
\[
\intrd\langle
\nabla_x\Phi(y',\omega)-\omega'\rangle^{-N}\,d\omega\leq
C',\quad \forall
(y',\omega')\in\R^d\times\R^d.
\]
To this end, we perform the change of
variable
$\R^d\ni\omega\longmapsto\nabla_x\Phi(y',\omega)$
which is a global diffeomorphims of
$\R^d$ in view \eqref{crescita2bis} and
\eqref{i1tris}. Moreover the Jacobian
determinant of its inverse is uniformly
bounded with respect to $y'$ (see the
discussion at the beginning of the
present section). Hence, the last
integral is, for $N>d$,
\begin{align*}
&\lesssim
 \int_{\R^{d}}
\langle\tilde{\omega}-\omega'\rangle^{-N}
\,d\tilde{\omega}=C'.
\end{align*}
The proof of \eqref{schur2} is
analogous and left to the reader.\par
Finally, the uniformity of the norm of
$A$ as a bounded operator, established
in the last part of the statement,
follows from the proof itself.
\end{proof}

\section{Sufficient Conditions for Boundedness of \psdo s}\label{sectionb}
Here we study the boundedness on
modulation spaces of pseusodifferential
operators, namely operators of the form
\eqref{pseudoformula} above, for some
symbol classes. First we consider the
case of symbols in $\sg^{m_1,m_2}$.
\par We observe that the full
pseudodifferential calculus is
available for these operators. Indeed,
it is a special case of the calculus
for general H\"ormander's classes
$S(m,g)$ associated with a weight $m$
and a metric $g$ (\cite[Chapter
XVIII]{hormander}). Here
$m(x,\eta)=\langle
x\rangle^{m_2}\langle\eta\rangle
^{m_1}$ and
$g_{x,\eta}(z,\zeta)=\langle
x\rangle^{-2}|dz|^2+\langle
\eta\rangle^{-2}|d\zeta|^2$. In
particular the composition of two
pseudodifferential operators with
symbols in $\sg^{m_1,m_2}$ and
$\sg^{m'_1,m'_2}$ is a
pseudodifferential operators with
symbol in $\sg^{m_1+m'_1,m_2+m'_2}$.\\
Now, we claim that such an operator,
with symbol in $\sg^{m_1,m_2}$, extends
to a bounded operator
$\tilde{M}^{p,q}_{s_1,s_2}\to
\tilde{M}^{p,q}_{s_1-m_1,s_2-m_2}$, for
every $s_1,s_2,m_1,m_2\in\R$, $1\leq
p,q\leq\infty$. This was proved in
\cite[Corollary 4.7]{Toftweight} when
$m_1=m_2=0$. When, in addition,
$s_1=s_2=0$ (the unweighted case) this
result is contained in \cite[Theorem
14.5.2]{book}. Our claim follows from
this latter special case by arguing as
follows.\par First observe that, for
 every $s_1,s_2\in\R$, the $\sg$
  pseudodifferential operator
 $\Lambda_{s_1,s_2}=\langle x\rangle^{s_2}\langle D\rangle^{s_1}
 $, of order $(s_1,s_2)$, is bounded (in fact it defines
 an isomorphism) from $\tilde{\mathcal M}^p_{s_1,s_2}$
 to $\tilde{\mathcal M}^p$ \cite[Theorem
 2.4]{Toftweight}. Moreover,
 $\Lambda_{s_1,s_2}^{-1}u=
 \langle D\rangle^{-s_1}(\langle x\rangle^{-s_2} u)$
 is a $\sg$ pseudodifferential
 operator of order
 $(-s_1,-s_2)$. By the above quoted composition formula,
   $$A=\Lambda_{m_1-s_1,m_2-s_2}
   A'\Lambda_{s_1,s_2},$$
    for a suitable $\sg$ pseudodifferential operator
   $A'$ of order $(0,0)$. As we already
   observed, $A'$ is bounded
   $\tilde{M}^{p,q}\to\tilde{M}^{p,q}$ by \cite[Theorem
14.5.2]{book},
   which gives the claim.\par
We now will show a more general
continuity result, for rougher symbols
on $\rdd$ satisfying estimates of the
type
\begin{equation}\label{simboli}
|\partial^\alpha_\eta\partial^\beta_x
\sigma(x,\eta)|\leq C_{\alpha,\beta}
\langle\eta\rangle^{m_1}\langle
x\rangle^{m_2},\quad |\alpha|\leq 2
N_2,\ |\beta|\leq 2 N_1,
\end{equation}
with
$\partial^\alpha_\eta\partial^\beta_x$
standing for distributional
derivatives. For $s_1,s_2\geq 0$, we
recall the definition
$v_{s_1,s_2}(x,\eta)=\la x\ra^{s_2}\la
\o\ra^{s_1}$.
 Our result reads as follows.
 \begin{theorem}\label{contmp}
For $s_1,s_2\geq 0$, let
$\mu\in\cM_{v_{s_1,s_2}}$.\par $(a)$
Consider a symbol $\sigma$ satisfying
\eqref{simboli}, with
$N_1>(d+s_1+|m_1|)/2$, $N_2>(d+s_2)/2$.
Then, for every $1\leq p,q\leq\infty$,
$\sigma(x,D)$ extends to a continuous
operator from $\tilde{M}^{p,q}_{\mu}$
to
$\tilde{M}^{p,q}_{{\mu}{v_{-m_1,-m_2}}}$.
\par
$(b)$ Consider a symbol $\sigma$
satisfying \eqref{simboli}, with
$N_1>(d+s_1)/2$, $N_2>(d+s_2+|m_2|)/2$.
Then, for every $1\leq p,q\leq\infty$,
$\sigma(x,D)$ extends to a continuous
operator from $\tilde{M}^{p,q}_{\mu
v_{m_1,m_2}}$ to
$\tilde{M}^{p,q}_{{\mu}{}}$.
\end{theorem}

To chase our goal, we first show an
approximate diagonalization of
$\sigma(x,D)$ by Gabor frames. In the
sequel, we consider a Gabor frame
$\{g_{k,n}\}_{k,n}$, $(k,n)\in
\a\zd\times\beta\zd$, with window
$g\in\mathcal{S}(\R^d)$. A small
variant of \cite[Theorem 3.1, Remark
3.2]{fio1} (see also
\cite{roch-tachiz98}) provides the
following almost diagonalization.
\begin{theorem}\label{diag} Consider a symbol
$\sigma$ satisfying \eqref{simboli}.
Then there exists $C_{N_1,N_2}>0$ such
that
\begin{equation}\label{ET1}
    |\la \sigma(x,D) g_{k,n},\, g_{k',n'}\ra|
    \leq C_{N_1,N_2}\frac{ \la n\ra^{m_1}\la k'\ra^{m_2}} {\la n-n'\ra^{2N_1}\la k-k'\ra^{2N_2}}.
\end{equation}
\end{theorem}
\begin{proof} The proof essentially
repeats that of \cite[Theorem 3.1,
Remark 3.2]{fio1}. Since that results
was actually established for more
general classes of FIOs, for the
convenience of the reader outline the
main ideas. An explicit computation
shows that
\begin{align*} |\langle \sigma(x,D)& (M_n T_k g),
M_{n'}
T_{k'}\gamma\rangle|\\
&=|\iint e^{2\pi i x(n-n')-\eta(k-k')}
\big[e^{2\pi i x\eta}
\sigma(x+k',\eta+n)\big]
\bar{\gamma}(x)\hat{g}(\o)\,dx
   d\o|.
   \end{align*}

Then one uses the identity
\begin{align*}
(1-\Delta_x)^{N_1} (1-\Delta_\o)^{N_2}&
 e^{2\pi i [x(n-n')-\o(k-k')]}\\
 &=\la
    2\pi(k-k')\ra^{2N_2}\!\la
    2\pi(n-n')\ra^{2N_1} \!e^{2\pi i [x(n-n')-\o(k-k')]},
\end{align*}
and integrates by parts. Since
$g\in\cS$, the estimates
\eqref{simboli} combined with Petree's
inequality $\langle
z+w\rangle^s\leq\langle
z\rangle^s\langle w\rangle^{|s|}$ give
 \eqref{ET1}.
\end{proof}\\
The proof of the boundedness property
of $\sigma(x,D)$
 makes use of the following generalization
 of the
  Schur Test, contained in \cite[Proposition
  5.1]{fio1}.
\begin{proposition}\label{proschur}
Consider an operator defined on
sequences on the lattice
$\Lambda=\alpha\mathbb{Z}^d\times\beta\mathbb{Z}^d$
by
\[
(Kc)_{m',n'}=\sum_{m,n}
{K_{m',n',m,n}}c_{m,n}.
\]
Assume
$$\{K_{m',n',m,n}\}\in \ell^\infty_{n}\ell^1_{n'}\ell^\infty_{m'}\ell^1_{m}\cap
\ell^\infty_{n'}\ell^1_{n}\ell^\infty_{m}\ell^1_{m'}\quad
\mbox{and}\quad \{K_{m',n',m,n}\}\in
\ell^\infty_{m',n'}\ell^1_{m,n}\cap
\ell^\infty_{m,n}\ell^1_{m',n'}.$$
Then, for every $1\leq p,q\leq\infty$,
the operator $K$ is continuous on
$\ell^{p,q}$ and on
$\tilde{\ell}^{p,q}$ (recall that
$\tilde{\ell}^{p,q}$ is the closure of
the space of eventually zero sequences
in ${\ell}^{p,q}$).
\end{proposition}
\vskip0.3truecm \noindent \textsl{Proof
of Theorem \ref{contmp}}. $(a)$
Consider a normalized tight frame
$\mathcal{G}(g,\a,\beta)$ with
$g\in\cS(\rd)$. In view of Theorem
\ref{teomod2}, showing the boundedness
of $\sigma(x,D)$ from
$\tilde{M}^{p,q}_{\mu}$ to
$\tilde{M}^{p,q}_{\mu v_{-m_1,-m_2}}$
is equivalent to proving the
boundedness of the infinite matrix
$$K_{k',n',k,n}=\la \sigma(x,D) g_{k,n},\, g_{k',n'}\ra \frac{\mu(k',n')}{\la k'\ra^{m_2}\la n'\ra^{m_1}\mu(k,n)}
$$
from $\tilde{\ell}^{p,q}$ into itself.
We make use of the Schur Test above
(Proposition \ref{proschur}). The
estimate \eqref{ET1} and the assumption
$\mu\in\cM_{v_{s_1,s_2}}$ combined with
Petree's inequality yield
$$|K_{k',n',k,n}|\lesssim \la n-n'\ra^{s_1+|m_1|-2N_1}\la
 k-k'\ra^{s_2-2N_2},
$$
so that
$$\sup_{k',n'\in\zd}\sum_{k,n\in\zd}|K_{k',n',k,n}|<\infty,
$$
because of the choice of $N_1,N_2$.
Analogously, one obtains
 $\{K_{m',n',m,n}\}\in\ell^\infty_{k,n}\ell^1_{k',n'}$. Similarly one obtains the estimate
\begin{align*}\sup_{n\in
\zd}\sum_{n'\in \zd}& \sup_{k'\in
\zd}\sum_{k\in
\zd}|K_{k',n',k,n}|\\
&\lesssim \sup_{k'\in
\zd}\sum_{k\in\zd}\la
k-k'\ra^{s_2-2N_2}\sup_{n\in
\zd}\sum_{n'\in\zd}\la
n-n'\ra^{s_1+|m_1|-2N_1}<\infty,
\end{align*}
that is,
$\{K_{m',n',m,n}\}\in\ell^\infty_n\ell^1_{n'}\ell^\infty_{k'}\ell^1_k$,
and also
$\{K_{m',n',m,n}\}\in\ell^\infty_{n'}\ell^1_{n}\ell^\infty_{k}\ell^1_{k'}$,
as desired.\par The proof of part $(b)$
is very similar and left to the reader.
\endproof
\begin{remark}\rm The formula
\eqref{ET1} is not symmetric with
respect to the space variables and the
dual variables. This is due to the fact
that we are using the so called ``left"
or Kohn-Nirenberg quantization
\eqref{pseudoformula}. Instead, the
Weyl quantization
(\cite{folland89,hormander}):
\[
\sigma^w(x,D)u= \iint e^{2\pi
i(x-y)\eta}
\sigma\left(\frac{x+y}{2},\eta\right)f(y)\,dy\,d\eta,
\]
combined with properties of the
cross-Wigner distribution as in
\cite{roch-tachiz98}, yields
\[
 | \la \sigma^w(x,D) g_{k,n}, g_{k',n'}\ra|\leq
 C_{N_1,N_2}
 \frac{\left\la{n+n'}\right\ra^{m_1}
 \left\la{k+k'}\right\ra^{m_2}}{\la n-n'\ra^{2N_1}
 \la k-k'\ra^{2N_2}}.
 \]
 Theorem \ref{contmp} for Weyl operators reads as follows:
\textsl{Let $s_1,s_2\geq 0$,
$\mu\in\cM_{v_{s_1,s_2}}$, $\sigma$ be
a symbol satisfying \eqref{simboli},
with $N_i>(d+s_i+|m_i|)/2$, $i=1,2$.
Then, for every $1\leq p,q\leq\infty$,
$\sigma(x,D)$ extends to a continuous
operator from $\tilde{M}^{p,q}_{\mu}$
to
$\tilde{M}^{p,q}_{{\mu}{v_{-m_1,-m_2}}}$
and from $\tilde{M}^{p,q}_{\mu
v_{m_1,m_2}}$ to
$\tilde{M}^{p,q}_{{\mu}{}}$.}
\end{remark}

\section{Proof of Theorem \ref{maintheorem}}\label{section5}
 The claim of Theorem
\ref{maintheorem} will follow if we
prove the boundedness of every $\sg$
FIO $A$ of order
$(m_1,m_2)=(-d/2,-d/2)$ on the endpoint
cases $M^1$ and
 $\tilde{M}^\infty$. Indeed, since it is known
from Theorem \ref{elle2} that a FIO of
order $(0,0)$ is bounded on
 $L^2=M^2$, the desired
 continuity result on $M^p$, when
 $m_1=m_2=-d|1/2-1/p|$,
 $1<p<\infty$,
  is due to complex
 interpolation, detailed below.\par
 As already observed, for
 every $s_1,s_2\in\R$, the operator
 $\Lambda_{s_1,s_2}u=\langle D\rangle^{s_1}
 (\langle x\rangle^{s_2} u)$ defines
 an isomorphism from $\tilde{M}^p_{s_1,s_2}$
 to $\tilde{M}^p$ \cite[Theorem 2.4]{Toftweight}. Its inverse
 $\Lambda_{s_1,s_2}^{-1}u=\langle x\rangle^{-s_2} (\langle D\rangle^{-s_1} u)$
 is a $\sg$ pseudodifferential
 operator of order $(-s_1,-s_2)$. If $A$ is a $\sg$ FIO of order
$(m_1,m_2)$, writing
\begin{equation}\label{fiodec}
A=T\Lambda_{s_1,s_2},\quad
T:=A\Lambda_{s_1,s_2}^{-1},
\end{equation}
by Theorem \ref{composition2}, the
operator $T$ is a FIO with the same
phase as the operator $A$ and of order
$(m_1-s_1,m_2-s_2)$. Now, assume that
Theorem \ref{maintheorem} is true for
$p=1,2$. Consider $1<p<2$ and let $A$
be a FIO of order $(m_1,m_2)$, with
$m_1=m_2=-d(1/p-1/2)$. For $p=1$,
provided that $s_i=m_i+d/2$, $i=1,2$,
by \eqref{fiodec} the operator $T$ is
of order $(-d/2,-d/2)$, hence bounded
on $M^1$. As a consequence,
  the FIO $A$ extends to a
 bounded operator from $M^1_{m_1+d/2,m_2+d/2}$ to $M^1$. For $p=2$ and $s_i=m_i$, $i=1,2$, the FIO $T$ is of order $(0,0)$, so that $T$ is bounded on $L^2$ and $A$ is bounded from $M^2_{m_1,m_2}$ to $M^2$. \\
 By complex interpolation (see \eqref{interpmp}), the operator $A$ is bounded between the following spaces $$ A: M^p_{\tilde{m}_1,\tilde{m}_2}=(M^1_{m_1+d/2,m_2+d/2},M^2_{m_1,m_2})_{[\theta]}\to M^p=(M^1,M^2)_{[\theta]},$$ with
 $1/p=(1-\theta)/1+\theta/2$, $\theta\in(0,1)$, $\tilde{m_i}=(1-\theta)(m_i+d/2)+\theta m_i=m_i+(1-\theta)d/2$, $i=1,2$.
 These equalities yield $\tilde{m_i}=m_i+d(1/p-1/2)=0$, because $m_i=-d(1/p-1/2)$ by assumption.
 Hence $A$ is bounded from $M^p$ to $M^p$, as desired. The
 proof for $2<p<\infty$ is
 similar.\par Of course, when
 in one of the inequalities in
 \eqref{soglia} (or in both)
  a strict inequality holds, the desired result
  follows from the
  equality-case, for an operator
  of order $(m'_1,m'_2)$
  with $m'_1\leq m_1$, $m'_2\leq m_2$, has also order $(m_1,m_2)$. \par
Hence, from now on, we assume
$m_1=m_2=-d/2$ and prove the
boundedness of $A$ on $M^1$ and on
$\tilde{M}^\infty$.

\subsection{Boundedness on
$M^1$}

 Consider now the usual
Littlewood-Paley decomposition of the
frequency domain. Namely,
 fix a smooth function $\psi_0(\eta)$
  such that $\psi_0(\eta)=1$
  for $|\eta|\leq1$ and
  $\psi_0(\eta)=0$ for
  $|\eta|\geq2$. Set
  $\psi(\eta)=\psi_0(\eta)-\psi_0(2\eta)$,
  $\psi_j(\eta)=\psi(2^{-j}\eta)$, $j\geq1$.
  Then
  \[
  1=\sum_{j=0}^\infty\psi_j(\eta),\quad
\forall\eta\in\R^d.
 \]
 Following the general philosophy
 of \cite{fefferman}, we perform a dyadic
 decomposition of the symbol
 $\sigma$ on boxes of size $2^k\times
 2^j$, $k,j\geq0$,
 hence tailored to the $\sg$ symbol
 estimates; then
  we conjugate each dyadic operator with dilations
   in such a way to transform any box
   into a cube of
   size $2^{(j+k)/2}\times2^{(j+k)/2}$.
   Finally we will apply
   Proposition \ref{pro2} to these transformed operators. \par
   Namely,
   consider the decomposition

\begin{equation}\label{somma}
A=\sum_{j,k\geq0} \tjk=\sum_{0\leq j<k}
\tjk+\sum_{0\leq k\leq j} \tjk,
\end{equation}
where $\tjk$ is the FIO with the same
phase $\Phi$ as $A$ and symbol
\begin{equation}\label{sommabis}
\sigma_{j,k}(x,\eta):=\psi_k(x)\sigma(x,\eta)
\psi_j(\eta).
\end{equation}
The key point here is that the symbols
$\sigma_{j,k}$ are supported where
$\langle\eta\rangle\asymp 2^j$,
$\langle x\rangle\asymp 2^{k}$ and, for
any $\a,\beta\in\zd_+$,
$$|\partial^\a_\o\partial^\beta_x \sigma_{j,k}\phas|\leq C_{\a,\beta}\langle \o\rangle^{-(\frac{d}2+|\a|)}\langle x\rangle^{-(\frac{d}2+|\beta|)},\quad \forall j,k\geq 0,
$$
i.e., $\sigma_{j,k}$ lie in a
bounded subset of $\sg^{-d/2,-d/2}$.\\
Moreover, we observe that
\begin{equation}\label{coniugazione}
\tjk=U_{2^{\frac{j-k}{2}}}\ttjk
U_{2^{-\frac{j-k}{2}}},
\end{equation}
where $\ttjk$ is the FIO with phase
\begin{equation}\label{i100bis}
{\Phi}_{j,k}(x,\eta):=\Phi(2^{-\frac{j-k}{2}}x,2^{\frac{j-k}{2}}
\eta),
\end{equation}
and symbol
\[
\tilde{\sigma}_{j,k}(x,\eta):=
\sigma_{j,k}(2^{-\frac{j-k}{2}}x,2^{\frac{j-k}{2}}\eta),
\]
and $U_\lambda f(y)=f(\lambda y)$,
$\lambda>0$, is the dilation
operator.\par Notice that
$\tilde{\sigma}_{j,k}$ is supported in
the set
\begin{equation}\label{nuc}
\mathcal{V}_C=\{(x,\eta)\in\R^{2d}:
C^{-1}2^j\leq \langle
2^{\frac{j-k}{2}}\eta\rangle \leq
C2^j,\ C^{-1}2^k\leq \langle
2^{\frac{k-j}{2}} x\rangle \leq C2^k\},
\end{equation}
for some $C>0$.\par We first consider
the sum over $k\leq j$ in
\eqref{somma}.\par Assume for a moment
that the following estimate holds
\begin{equation}\label{intermedia}
\|\ttjk u\|_{M^1}\lesssim
2^{-(j+k)d/2}\|u\|_{M^1},
\end{equation}
 and recall the dilation properties for modulation spaces (Theorem \ref{dilprop}), for $p=1$:
\begin{equation}\label{di1}
\|U_{\lambda}f\|_{M^1}\lesssim
\|f\|_{M^1},\quad \lambda\geq1,
\end{equation}
and
\begin{equation}\label{di2}
\|U_{\lambda}f\|_{M^1}\lesssim
\lambda^{-d}\|f\|_{M^1},\quad
0<\lambda\leq1.
\end{equation}
Then, using \eqref{di1} (with
$\lambda=2^{(j-k)/2}$) and \eqref{di2}
(with $\lambda=2^{-(j-k)/2}$) we obtain
\[
\|\tjk u\|_{M^1}\lesssim
2^{-kd}\|u\|_{M^1}.
\]
Actually, for the frequency
localization of $\tjk$, the following
finer estimate holds:
\begin{equation}\label{boi}
\|\tjk u\|_{M^1}=\|\tjk(\chi(2^{-j}D)
u)\|_{M^1}\lesssim
2^{-kd}\|\chi(2^{-j}D)u\|_{M^1},\
j\geq1,
\end{equation}
where $\chi$ is a smooth function
satisfying $\chi(\eta)=1$ for $1/2\leq
|\eta|\leq2$ and $\chi(\eta)=0$ for
$|\eta|\leq1/4$ and $|\eta|\geq4$ (so
that $\chi\psi=\psi$). Summing on $j,k$
this last estimate, with the aid of
Lemma \ref{le41} $(a)$, we obtain
\begin{align*}
\|\sum_{0\leq k\leq j} \tjk
u\|_{M^1}&\leq\sum_{j=0}^\infty\sum_{k=0}^j
\|\tjk u\|_{M^1}\\
&\lesssim
\|u\|_{M^1}+\sum_{j=1}^\infty\sum_{k=0}^j
 2^{-kd}
 \|\chi(2^{-j}D)u\|_{M^1}\\
 & \lesssim\|u\|_{M^1}+\sum_{j=1}^\infty
 \|\chi(2^{-j}D)u\|_{M^1}\\
 &\lesssim \|u\|_{M^1}.
\end{align*}
which is the desired estimate for the
sum over $k\leq j$.\par It remains to
prove \eqref{intermedia}. This follows
from Proposition \ref{pro2} applied to
the operator $2^{(j+k)d/2}\ttjk$.
Indeed, it is easy to see that the
hypotheses are satisfied uniformly with
respect to $j,k$. Precisely, the chain
rule gives, for every $j,k\geq 0$,
\[
|\partial^\alpha_\eta\partial^\beta_x
\tilde{\sigma}_{j,k}(x,\eta)|\lesssim
2^{(j+k)\left(-\frac{d}{2}-\frac{|\alpha|+|\beta|}{2}
\right)}.
\]
Similarly, we also have
\begin{equation}\label{i10bis}
|\partial^\alpha_\eta\partial^\beta_x
{\Phi}_{j,k}(x,\eta)|\lesssim
2^{(j+k)\left(1-\frac{|\alpha|+|\beta|}{2}\right)},\quad
{\rm for\ every}\
(x,\eta)\in\mathcal{V}_{C'},
\end{equation}
for every fixed $C'>0$ (see \eqref{nuc}
and notice that $\mathcal{V}_{C'}$
contains an $\epsilon$-neighborhood of
the support of $\tilde{\sigma}_{j,k}$
if $C'$ is large and $\epsilon$ small
enough). Clearly we also have
\begin{equation}\label{crescita2tris}
|\partial^\alpha_\eta\partial^\beta_x
\Phi_{j,k}(x,\eta)|\leq
C_{\alpha,\beta},\quad
|\alpha|=|\beta|=1,\
\forall(x,\eta)\in\R^{2d},
\end{equation}
 and
\begin{equation}\label{i1bis}
\left|{\rm det}\,
\left(\frac{\partial^2\Phi_{j,k}}{\partial
x_i\partial \eta_l}\Big|_{
(x,\eta)}\right)\right|>\delta>0,\quad
\forall (x,\eta)\in\R^{2d}.
\end{equation}
Hence Proposition \ref{pro2}
 applies and gives, for $1\leq
 p\leq\infty$,
 \begin{equation}\label{intermedia2}
\|\ttjk u\|_{M^p}\lesssim
2^{-(j+k)d/2}\|u\|_{M^p}.
\end{equation}
For $p=1$ this is
 \eqref{intermedia}.\par
 We now consider the sum over
 $j< k$ in
 \eqref{somma}. Namely, we prove that
 \begin{equation}
 \|\sum_{0\leq j<k}\tjk
 u\|_{M^1}\lesssim\|u\|_{M^1}.
\end{equation}
Using the Littlewood-Paley
decomposition $\sum_{l\geq
0}\psi_l(x)=1$, we write
\[
\sum_{0\leq j<k}\tjk
 u=\sum_{l=0}^\infty \sum_{0\leq j<k}\tjk
 (\psi_l u).
 \]
By the triangle inequality and Lemma
 \ref{le41} $(b)$, it
 suffices to prove
 \begin{equation}\label{lal}
\|\sum_{0\leq j<k}\tjk
 (\psi_l u)\|_{M^1}\lesssim
 \|u\|_{M^1}.
 \end{equation}
More precisely, one should apply this
estimate with $u$ replaced by
$\chi(2^{-l}\cdot)u$, with $\chi$ as in
\eqref{boi} above and then use Lemma
\ref{le41} $(b)$.\par Applying Theorem
\ref{composition2} with
$p(x,D)=\psi_l(x)$ (notice that the
multiplicative operators $\psi_l(x)$
are pseudodifferential operators with
symbols in a bounded subset of
$\sg^{(0,0)}$) to each composition
 $2^{(j+k)d/2}\tjk \psi_l$, we obtain
 \[
\tjk\psi_l=2^{-(j+k)d/2}S_{j,k}^{(l)}
+2^{-(j+k)d/2}R_{j,k}^{(l)}
\]
where $S_{j,k}^{(l)}$ are FIOs with the
same phase $\Phi_{j,k}$ in
\eqref{i100bis} and symbols
$\sigma_{j,k}^{(l)}$ belonging to
bounded subset of $\sg^{0,0}$,
supported in
\begin{equation}\label{b1bis}
\{(x,\eta)\in\R^{2d}: \langle
\eta\rangle\asymp 2^j,\
\langle\nabla_\eta\Phi(x,\eta)\rangle\asymp
2^l, \langle x\rangle\asymp 2^{k}\}.
\end{equation}
The operators $R_{j,k}^{l}$ are
smoothing operators whose symbols
$r_{j,k,l}$ lie
 in a bounded subset of
 $\mathcal{S}(\R^{2d})$.
 \par Observe that,
by \eqref{crescita},
\[
\langle\nabla_\eta\Phi(x,\eta)\rangle\asymp\langle
x \rangle.
\]
Inserting this equivalence in
\eqref{b1bis}, we deduce that there
exists $N_0>0$ such that
$\sigma_{j,k}^{(l)}$ vanishes
identically if $|k-l|>N_0$. Whence, the
left-hand side in \eqref{lal} is seen
to be
\[
\leq\sum_{k\geq0: |k-l|\leq
N_0}\sum_{j=0}^{k-1}
2^{-(j+k)d/2}\|S_{j,k}^{(l)}u \|_{M^1}+
\sum_{k=0}^\infty\sum_{j=0}^{k-1}
2^{-(j+k)d/2}
\|R_{j,k}^{(l)}u\|_{M^1}.\] Since
\[
\|R_{j,k}^{(l)}u\|_{M^1}\lesssim
\|u\|_{M^1},
\]
\eqref{lal} will follow from
\begin{equation}\label{parteprincipale2}
\|S_{j,k}^{(l)}u \|_{M^1}\lesssim
2^{kd/2} \|u\|_{M^1}.
\end{equation}
In order to prove this estimate, we
write
\[
S_{j,k}^{(l)}:=U_{2^{\frac{j-k}{2}}}
\tilde{S}_{j,k}^{(l)}
U_{2^{-\frac{j-k}{2}}},
\]
where $\tilde{S}_{j,k}^{(l)}$ is the
FIO with phase ${\Phi}_{j,k}(x,\eta)$
defined in \eqref{i100bis}, and symbol
\[
\tilde{\sigma}_{j,k}^{(l)}(x,\eta):=
\sigma_{j,k}^{(l)}(2^{-\frac{j-k}{2}}x,2^{\frac{j-k}{2}}
\eta),
\]
supported in a set $\mathcal{V}_C$ of
the type \eqref{nuc}.\\ Now, taking
into account \eqref{di1}, \eqref{di2},
we see that \eqref{parteprincipale2}
will follow (with an additional factor
$2^{-jd/2}$) from
\begin{equation}\label{lastebis}
\|\tilde{S}_{j,k}^{(l)}
u\|_{M^1}\lesssim \|u\|_{M^1}.
\end{equation}
Precisely, if \eqref{lastebis} holds
true, we have
\begin{equation*}
\|S_{j,k}^{(l)}u \|_{M^1}\lesssim
2^{-d\frac{j-k}{2}}\|\tilde{S}_{j,k}^{(l)}U_{2^{-\frac{j-k}{2}}}u\|_{M^1}
\lesssim
2^{-d\frac{j-k}{2}}\|U_{2^{-\frac{j-k}{2}}}u\|_{M^1}\lesssim
2^{-d\frac{j-k}{2}}\|u\|_{M^1}.
\end{equation*}
The estimate \eqref{lastebis} is a
consequence of Proposition \ref{pro2}
applied to $\tilde{S}_{k,j}$. Indeed,
since the symbols $\sigma_{j,k}^{(l)}$
belong to a bounded subset of
$\sg^{0,0}$ and are supported where
$\langle\eta\rangle\asymp 2^j$,
$\langle x\rangle\asymp 2^{k}$, it
turns out
\[
|\partial^\alpha_\eta\partial^\beta_x
\tilde{\sigma}_{j,k,l}(x,\eta)|\lesssim
2^{-(j+k)\frac{|\alpha|+|\beta|}{2}}.
\]
 On the
other hand, we already observed that
\eqref{i10bis}, \eqref{crescita2tris}
and \eqref{i1bis} hold true. Hence
Proposition \ref{pro2} gives
\eqref{lastebis}.\par

 \subsection{Boundedness on
 $\tilde{M}^\infty$}

We now show the boundedness of $A$ (of
order $(m_1,m_2)=(-d/2,-d/2)$) on
$\tilde{M}^\infty$, using the notations
of the previous subsection. Our
arguments in this case reflect the
symmetry of the $\sg$ symbol classes
with respect to the exchange
$x\longleftrightarrow \o$.

Consider again the decomposition in
\eqref{somma}. Hence, $\tjk$ is the FIO
with phase $\Phi_{j,k}$ in
\eqref{i100bis} and symbol
$\sigma_{j,k}$ in \eqref{sommabis}. We
first teat the sum over $k\leq j$. By
Lemma \ref{lemma2} $(a)$ we have
\begin{align}
\|\sum_{0\leq k\leq j}\tjk
u\|_{M^\infty}&=\|\sum_{l\geq0}\psi_l(D)\sum_{0\leq
k\leq j}\tjk u\|_{M^\infty}
\nonumber\\
&\lesssim\sup_{l\geq0}\|\psi_l(D)\sum_{0\leq
k\leq j}\tjk u\|_{M^\infty}
\nonumber\\
&\leq\sup_{l\geq0}\sum_{0\leq k\leq
j}\|\psi_l(D)\tjk
u\|_{M^\infty}.\label{continuazione}
\end{align}
 Applying Theorem
\ref{composition} to each product
$\psi_l(D) 2^{(j+k)d/2}\tjk$, we have
\[
\psi_l(D)\tjk=
2^{-(j+k)d/2}S_{j,k,l}+2^{-(j+k)d/2}R_{j,k,l},
\]
where $S_{j,k,l}$ are FIOs with the
same phase $\Phi_{j,k}$ and symbols
$\sigma_{j,k,l}$ belonging to a bounded
subset of $\sg^{0,0}$, supported in
\begin{equation}\label{b1}
\{(x,\eta)\in\R^{2d}: \langle
x\rangle\asymp 2^k,\
\langle\nabla_x\Phi(x,\eta)\rangle\asymp
2^l, \langle\eta\rangle\asymp 2^{j}\}.
\end{equation}
The operators $R_{j,k,l}$ are smoothing
operators whose symbols $r_{j,k,l}$ are
 in a bounded subset of
 $\mathcal{S}(\R^{2d})$.
 \par Observe that,
by \eqref{crescita},
\[
\langle\nabla_x\Phi(x,\eta)\rangle\asymp\langle\eta\rangle.
\]
Inserting this equivalence in
\eqref{b1}, we obtain that there exists
$N_0>0$ such that $\sigma_{j,k,l}$
vanishes identically if $|j-l|>N_0$.
Whence, the right-hand side in
\eqref{continuazione} is seen to be
\[
\leq\sup_{l\geq0}\sum_{j\geq0:
|j-l|\leq N_0}\sum_{k=0}^j
2^{-(j+k)d/2}\|S_{j,k,l}u
\|_{M^\infty}+\sup_{l\geq0}
\sum_{j=0}^\infty\sum_{k=0}^j
2^{-(j+k)d/2}
\|R_{j,k,l}u\|_{M^\infty}.\] This
expression will be dominated by the
$M^\infty$ norm of $u$ if we prove that
\begin{equation}\label{parteprincipale}
\|S_{j,k,l}u\|_{M^\infty}\lesssim
2^{jd/2}\|u\|_{M^\infty},
\end{equation}
because clearly
\begin{equation}\label{resto}
\|R_{j,k,l}u\|_{M^\infty}\lesssim
\|u\|_{M^\infty}.
\end{equation}
 To prove
\eqref{parteprincipale} we recall from
Theorem \ref{dilprop} that
\begin{equation}\label{di1bis}
\|U_{\lambda}f\|_{M^\infty}\lesssim
\|f\|_{M^\infty},\quad \lambda\geq1,
\end{equation}
and
\begin{equation}\label{di2bis}
\|U_{\lambda}f\|_{M^\infty}\lesssim
\lambda^{-d}\|f\|_{M^\infty},\quad
0<\lambda\leq1.
\end{equation}
Then we write
\[
S_{j,k,l}:=U_{2^{\frac{j-k}{2}}}
\tilde{S}_{j,k,l}
U_{2^{-\frac{j-k}{2}}},
\]
where $\tilde{S}_{j,k,l}$ is the FIO
with phase ${\Phi}_{j,k}(x,\eta)$ in
\eqref{i100bis}, and symbol
\[
\tilde{\sigma}_{j,k,l}(x,\eta):=
\sigma_{j,k,l}(2^{-\frac{j-k}{2}}x,2^{\frac{j-k}{2}}
\eta),
\]
supported in a set $\mathcal{V}_C$ of
the type \eqref{nuc}.\\ Now, taking
into account \eqref{di1bis},
\eqref{di2bis}, we see that
\eqref{parteprincipale} will follow
(with an additional factor $2^{-kd/2}$)
from
\begin{equation}\label{laste}
\|\tilde{S}_{j,k,l}
u\|_{M^\infty}\lesssim
\|u\|_{M^\infty}.
\end{equation}
This last estimate is a consequence of
Proposition \ref{pro2} applied to
$\tilde{S}_{k,j}$. Indeed, since the
symbols $\sigma_{j,k,l}$ belong to a
bounded subset of $\sg^{0,0}$ and are
supported where
$\langle\eta\rangle\asymp 2^j$,
$\langle x\rangle\asymp 2^{k}$, it
turns out
\[
|\partial^\alpha_\eta\partial^\beta_x
\tilde{\sigma}_{j,k,l}(x,\eta)|\lesssim
2^{-(j+k)\frac{|\alpha|+|\beta|}{2}}.
\]
Again, \eqref{i10bis},
\eqref{crescita2tris} and \eqref{i1bis}
have already been verified. Hence
Proposition \ref{pro2} gives
\eqref{laste}.\par We now treat the sum
over $j<k$ in \eqref{somma}.\\
By Lemma \ref{lemma2} $(b)$ and the
triangle inequality we have
\begin{align*}
\|\sum_{0\leq j<k }\tjk
u\|_{M^\infty}&=\|\sum_{k=0}^\infty\sum_{j=0}^{k-1}\tjk
u\|_{M^\infty}\\
&\lesssim
\sup_{k\geq0}\sum_{j=0}^{k-1}\|\tjk
u\|_{M^\infty}.
\end{align*}
Hence, the desired result will follow
from the estimate
\[
\|\tjk u\|_{M^\infty}\lesssim
2^{-jd}\|u\|_{M^\infty}.
\]
Using \eqref{coniugazione} and
\eqref{di1bis}, \eqref{di2bis}, we see
that it suffices to prove
\[
\|\ttjk u\|_{M^\infty}\lesssim
2^{-\frac{(j+k)d}{2}}\|u\|_{M^\infty},
\]
but this is \eqref{intermedia2} for
$p=\infty$. \par This concludes the
proof of Theorem \ref{maintheorem}.
\begin{remark}\rm
Theorem \ref{maintheorem} holds true
for operators of Type II as well (see
\eqref{FIOB}), as one sees by using
\eqref{aggiunto}.
\end{remark}

\section{Sharpness of the results and negative results for $L^p$}\label{sharp}
In this section we prove the sharpness
of Theorems \ref{maintheorem}.
Precisely, if one of the index pairing
$m_1, m_2$ fulfills $m_i
>\displaystyle -d
\left|\frac12-\frac1p\right|$, $1\leq
p\leq\infty$, ($i=1,2$), there are FIOs
of the type \eqref{FIO} and order
$(m_1,m_2)$, satisfying the assumptions
in the Introduction, which do not
extend to bounded operators on $M^p$,
$1\leq p<\infty$, nor on
$\tilde{M}^\infty$.\par In fact in
\cite{cordero-nicola-rodino} we
exhibited, for every $1\leq
p\leq\infty$, $m>-d|1/2-1/p|$, a FIO
which does not extend to a bounded
operator on $\tilde{M}^p$, with the
following features. The phase
$\Phi(x,\eta)=\sum_{j=1}^d\f(x_j)\eta_j$
is linear in $\eta$, where $\f:\R\to\R$
is a diffeomorphism satisfying
\eqref{bound} below. The symbol
$\sigma(x,\eta)$ belongs to
H\"ormander's class $S^m_{1,0}$ and is
compactly supported with respect to
$x$. In particular we see that $\Phi\in
\sg^{1,1}$ and satisfies
\eqref{nondeg}, and
$\sigma\in\sg^{m,-\infty}$. This shows
that the threshold for the index $m_1$
in Theorem \ref{maintheorem} is sharp,
even for symbols compactly supported in
$x$. For the sake of completeness we
briefly recall the construction of such
an operator. Then we show that the
threshold for the index $m_2$ is sharp
as well, even for symbols which are
compactly supported with respect to
$\eta$. Finally we show how the example
in this latter case gives the following
negative result for $L^p$ spaces.
\begin{proposition}\label{casolp}
For every $1\leq p\leq\infty$,
$m>-d|1/2-1/p|$, there exists a FIO
having phase
$\Phi(x,\eta)=\sum_{k=1}^d\f(\eta_k)x_k$,
where $\f:\R\to\R$ is a diffeomorphism
satisfying \eqref{bound} below, and
symbol compactly supported with respect
to $\o$ and in the class
$\sg^{-\infty,m}$, which does not
extend to a bounded operator on $L^p$,
$1\leq p<\infty$, nor on the closure of
the Schwartz space in $L^\infty$, if
$p=\infty$.
\end{proposition}

We first recall some results of
\cite{cordero-nicola-rodino}.
\begin{proposition}\label{prop1} Let $\f: \R\to \R$ be a
$\cC^{\infty}$ diffeomorphism, whose
restriction to the interval $(0,1)$ is
a non-linear
 diffeomorphism on $(0,1)$. This means that there
  exists a point $t_0\in(0,1)$ such that
  $\f^{''}(t_0)\not=0$.
  Let $\chi\in\cC^{\infty}_0(\R)$, $\chi\geq 0$,
   with $\chi(\f(t_0))\not=0$. Then, if we set
\begin{equation}\label{fn}
f_n (t)= \chi(t) e^{2\pi i nt},\quad
n\in\N,
\end{equation}
for $1\leq p\leq 2$, we have
\begin{equation}\label{fnest}
\|f_n\circ\f\|_{\cF L^p} \geq c \,
n^{1/p-1/2},\quad \quad\forall n\in\N.
\end{equation}
\end{proposition}

The generalization to dimension $d\geq
1$ reads as follows.
\begin{corollary} Let $\f$ be as in Proposition \ref{prop1}
and $f_n$ defined in \eqref{fn}. We
define
\begin{equation}\label{fntilde}
\tilde{f}_n(t_1,\dots,t_d)=f_n(t_1)\cdots
f_n(t_d),\quad
\tilde{\f}(t_1,\dots,t_d)=(\f(t_1),\dots,\f(t_d)),
\end{equation}
then
\begin{equation}\label{estfn}
\| \tilde{f}_n\circ \tilde{\f}\|_{\cF
L^p(\rd)} \geq c \,n^{d(1/p-1/2)},
\end{equation}
for $1\leq p\leq 2$.
\end{corollary}
The action of the multiplier $\la
D\ra^m$ on the functions $\tilde{f}_n$
is the following.
\begin{lemma}
Let $m\in\R$ and $\tilde{f}_n$ defined
in \eqref{fntilde}. Then,
\begin{equation}\label{dfnf}
    \|\la D\ra^{m} \tilde{f}_n\|_{M^p}\leq C n^m.
\end{equation}
\end{lemma}
\medskip
We can now prove the sharpness of
Theorem \ref{maintheorem}. The key idea
is that a $C^1$ change of variables
that leaves the $\cF L^p$ spaces
invariant must be affine (the so-called
Beurling-Helson Theorem
\cite{beurling53,lebedev94,kasso07}).\par\medskip
\noindent {\bf Sharpness of the
threshold for the frequency index
$m_1$} (see
\cite{cordero-nicola-rodino} for
details). We first study the case
$1\leq p\leq2$. Consider the FIO
$$T_{\tilde{\f}} f(x)= f\circ \tilde{\f}(x)=
\intrd e^{2\pi i \tilde{\f}(x) \eta}
\hat{f}(\eta)\,d\eta,$$ where
$\tilde{\f}$ is defined in
\eqref{fntilde}. We require that the
one-dimensional diffeomorphism $\f$
satisfies the assumptions of
Proposition \ref{prop1} and the
additional hypothesis
\begin{equation}\label{bound}\f(x)=x,\quad \mbox{for}
\,\,\,|x|\geq 1.
\end{equation}
Then, the phase $\Phi(x,\eta)=
\tilde{\f}(x) \eta$ fulfills
$\Phi\in\sg^{1,1}$ and is
non-degenerate. Notice that
$T_{\tilde{\f}}$ maps
$\mathcal{C}^\infty_0(\R^d)$ into
itself and ${\rm
supp}\,T_{\tilde{\f}}f\subset (0,1)^d$
if ${\rm supp}\,{f}\subset (0,1)^d$.
\par

Let $G\in\cC^{\infty}_0(\rd)$,
 $G\geq 0$ and $G\equiv1
$ on $[0,1]^d$. For $m_1\in\R$, the
symbol $a(x,\o)=G(x)\la\eta\ra^{m_1}$
satisfies $a\in \sg^{m_1,-\infty}$, and
the related FIO is given by
\begin{equation}\label{fioF}
Af(x)=\intrd e^{2\pi i \tilde{\f}(x)
\eta} G(x)\la\eta\ra^{m_1}
\hat{f}(\o)\,d\o =G(x)[(T_{\tilde{\f}}
\la D\ra^{m_1} )f](x).
\end{equation}
If $m_1\leq-d|1/2-1/p|$, Theorem
\ref{maintheorem} assures the
boundedness of $A$ on $M^p$. We now
show that this threshold is sharp for
$1\leq p\leq2$. Indeed, consider the
functions $\tilde{f}_n$ in
\eqref{fntilde}. They are supported in
$(0,1)^d$, so
$T_{\tilde{\f}}\tilde{f}_n$ are. Hence,
applying the estimate \eqref{estfn} and
Lemma \ref{lloc}, we obtain
\begin{align*}
n^{d(1/p-1/2)}&\lesssim \|
\tilde{f}_n\circ \tilde{\f}\|_{\cF
L^p(\rd)}=\|
T_{\tilde{\f}}\tilde{f}_n\|_{\cF
L^p(\rd)}
=\| GT_{\tilde{\f}}\tilde{f}_n\|_{\cF L^p(\rd)}\\&\asymp \| GT_{\tilde{\f}}\tilde{f}_n\|_{M^p(\rd)}=\| GT_{\tilde{\f}}\la D\ra^{m_1}\la D\ra^{-m_1}\tilde{f}_n\|_{M^p(\rd)}\\
&\lesssim \|F\|_{M^p\rightarrow M^p}
\|\la
D\ra^{-m_1}\tilde{f}_n\|_{M^p(\rd)}\lesssim
\|F\|_{M^p\rightarrow M^p} \,n^{-m_1},
\end{align*}
where the last inequality is due to
\eqref{dfnf}. For $n\to\infty$, we
obtain $-m_1 \geq d(1/p-1/2)$, i.e.,
\eqref{soglia}. \par

\vskip0.3truecm We now study the case
$2<p\leq\infty$.
  Observe that the adjoint operator
$T^*_{\tilde{\f}}$ of the above FIO
$T_{\tilde{\f}}$ is still a FIO given
by
$$T^*_{\tilde{\f}} f(x)=\frac1{|J_{\tilde{\f}}
(\tilde{\f}^{-1}(x))|}\intrd e^{2\pi i
\widetilde{\f}^{-1}(x)
\eta}f(\o)\,d\o,$$ with
$\widetilde{\f}^{-1}(x_1,\dots,x_d)=
(\f^{-1}(x_1),\dots, \f^{-1}(x_d))$ and
$|J_{\tilde{\f}}|$ the Jacobian of
$\tilde{\f}$. Its phase $\Phi\phas
=\widetilde{\f}^{-1}(x) \eta$ still
fulfills $\Phi\in\sg^{1,1}$ and the
standard assumptions.\par
 Now, let
$H\in \cC^{\infty}_0(\rd)$ , $H\geq 0$,
and $H(x)\equiv 1$ on supp ($G\circ
\f^{-1}$). For $m_1\in\R$, we define
the operator
\begin{equation}\label{ftilde}
\tilde{A}f(x)=H(x)[\la D\ra^{m_1}
T^*_{\tilde{\f}} (G f)](x).
\end{equation}
Using Theorem \ref{composition}, it is
easily seen that $\tilde{A}$ is a FIO
with symbol in $\sg^{m_1,-\infty}$ (the
symbol is compactly supported in the
$x$-variable). Its adjoint is given by
\begin{equation}\label{ftildestar}
\tilde{A}^*= G T_{\tilde{\f}}\la
D\ra^{m_1} H= A+ R,
\end{equation}
where $A$ is defined in \eqref{fioF}
and the remainder $R$ is given by
\[
R f(x)=G(x) [T_{\tilde{\f}}\la
D\ra^{m_1} ((H-1)f)](x).\] If we choose
a function
$\tilde{G}\in\cC^{\infty}_0(\rd)$ ,
$\tilde{G}\equiv 1$ on supp $G$ we can
write
\begin{align*} Rf&=
\tilde{G}(x)G(x) [T_{\tilde{\f}}\la
D\ra^{m_1}
((H-1)f)](x)\\&=\tilde{G}(x)T_{\tilde{\f}}[(G\circ
\tilde{\f}^{-1})\la D\ra^{m_1}
((H-1)f)](x).
\end{align*}
 By assumptions, supp $(G\circ \f^{-1})\cap$
supp $(H-1)=\emptyset$, so that the
pseudodifferential operator
$$f\longmapsto (G\circ \f^{-1})\la D\ra^{m_1} ((H-1)f)$$
is a regularizing operator
 (it immediately follows by the composition
  formula of pseudodifferential operators, see e.g.
   \cite[Theorem 18.1.8, Vol. III]{hormander}):
 this means that it maps $\cS'(\rd)$ into $\cS(\rd)$. The operator $T_{\tilde{\f}}$ is a
 smooth change of variables, so $\tilde{G}(x)T_{\tilde{\f}}$ maps
$\cS(\rd)$ into itself. To sum up, the
remainder operator $R$ maps $\cS'(\rd)$
into $\cS(\rd)$, hence it is bounded on
$M^{p}$. This means that $\tilde{A}^*$
is continuous on some $M^{p}$ iff $A$
is.

The operator $\tilde{A}$ is
 a FIO, with symbol in
$\sg^{m_1,-\infty}$ (compactly
supported in the $x$ variable). Hence
it is bounded on $M^{p}$ if
$m_1\leq-d|1/2-1/p|$ fulfills
\eqref{soglia}. We now show that this
threshold is sharp for $2<p<\infty$.
Indeed, if
 $\tilde{A}$ were
bounded on $M^{p}$, then its adjoint
$\tilde{A}^*$ would be bounded on
$(M^{p})'=M^{p'}$, with $1<p'<2$, and
 the same for $A$. But the former case gives
 the boundedness of
$A$ on $M^{p'}$ iff $-m_1\geq
d(1/p'-1/2)=d(1/2-1/p)$,
 that is the
desired threshold. For $p=\infty$, if
$\tilde{A}$ were bounded on
$\tilde{M}^\infty$, its adjoint
$\tilde{A}^*$ would be bounded on
$(\tilde{M}^\infty)'=M^{1}$ and the
former argument applies.\par\medskip
\noindent \textbf{Sharpness of the
threshold for
 the space index $m_2$.} The argument rely on the previous
 counterexample, combined
 with
   the Fourier invariance
  of $M^p$ and a duality
  trick.\par
 Consider, for $1\leq
 p\leq\infty$, $m>-d|1/2-1/p|$
 the type I FIO
 $A=A_{\Phi,\sigma}$
 constructed in the previous
 subsection. Hence,
 $\Phi(x,\eta)=\sum_{k=1}^d\f(x_k)\eta_k$,
for a diffeomorphism $\f:\R\to\R$
satisfying
 \eqref{bound},
 $\sigma\in\sg^{m,-\infty}$ is
 compactly supported with
 respect to $x$, and $A$ does
 not extends to a bounded
 operator on $\tilde{M}^p$.\\
   Let us set
$$ {}^t\Phi\phas=\Phi(\o,x)\quad \sigma^*\phas=
\overline{\sigma(\o,x)}.$$ Then, by
comparing the two definitions
\eqref{FIO} and \eqref{FIOB}, we have
\begin{equation}\label{link}
B_{- {}^t\Phi,\sigma^*}=\cF\circ
A_{\Phi,\sigma} \circ \cF^{-1} ,
\end{equation}
where $B_{- {}^t\Phi,\sigma^*}$ is the
type II operator in \eqref{FIOB} having
phase $- {}^t\Phi$ and symbol
$\sigma^*$. Using \eqref{link} and the
fact that the Fourier transform defines
an isomorphism of any $\tilde{M}^p$ we
see that the operator $B_{-
{}^t\Phi,\sigma^*}$ does not extents to
a bounded operator on $\tilde{M}^p$.
The same therefore holds for $(B_{-
{}^t\Phi,\sigma^*})^\ast$ on
$\tilde{M}^{p'}$, since
$\tilde{M}^{p'}=(\tilde{M}^{p})'$. On
the other hand, by \eqref{aggiunto} we
have $(B_{-
{}^t\Phi,\sigma^*})^\ast=A_{-
{}^t\Phi,\sigma^*}$. The last operator
possesses symbol
$\sigma^*\in\sg^{-\infty,m}$, compactly
supported with respect to $\eta$, and
gives the desired
counterexample.\par\medskip\noindent
{\bf Proof of Proposition
\ref{casolp}}. We start with an
elementary remark. Consider a FIO $A$
and suppose that it does {\it not}
satisfy an estimate of the type
\[
\|Au\|_{M^p}\leq C\|u\|_{M^p},\quad
\forall u\in\mathcal{S}(\R^d).
\]
Suppose, in addition, that the
distribution kernel of $K(x,y)$ of $T$
has the property that the two
projections of ${\rm supp}\, K$ on
$\R^d_x$ and $\R^d_y$ are bounded sets.
Then, it follows by Lemma \ref{lloc}
$(i)$ that $A$ does not extend to a
bounded operator on $\Fur L^p$, if
$1\leq p<\infty$, nor on the closure of
the Schwartz functions in $\Fur
L^\infty$, if $p=\infty$.\par Taking
this fact into account, we see that the
operator $\tilde{A}$ in \eqref{ftilde},
if $m_1>-d\left|{1/2}-{1/p}\right|$,
does not extend to a bounded operator
on $\Fur L^{p'}$, $2<p'<\infty$, nor on
the closure of the Schwartz space in
$\Fur L^\infty$, if $p'=\infty$. This
operator has a phase
$\Phi(x,\eta)=\sum_{k=1}^d\f(x_k)\eta_k$,
for a diffeomorphism $\f:\R\to\R$
satisfying
 \eqref{bound}, and a symbol
 $\tau\in\sg^{m_1,-\infty}$,
 compactly supported with
 respect to $x$. By repeating
 the same arguments as in the
 proof of the sharpness of
 the space index $m_2$ we see
 that the operator $A_{-
{}^t\Phi,\tau^*}$, with phase
${}^t\Phi\phas=\Phi(\o,x)$ and symbol
$\tau^*\phas= \overline{\tau(\o,x)}$,
does not extend to a bounded operator
on $L^p$, $1\leq p<2$, and gives the
desired counterexample in that
case.\par Similarly, for $2\leq
p\leq\infty$, we consider the operator
$\tilde{A}^\ast$ in \eqref{ftildestar}
which, for the same reason, does not
extend to a bounded operator on $\Fur
L^{p'}$, $1\leq p'\leq2$, if
$m_1>-d|1/2-1/p|$. The same holds for
$A$ in \eqref{fioF}, because of the
second inequality in
\eqref{ftildestar}. By arguing as
before, we obtain the requested
counterexample, having the desired
phase and symbol $G(\eta)\langle
x\rangle^{m_1}$.\par
 Observe that this latter
  example is precisely that observed
  in the Introduction (see \eqref{newsy}). Actually, in Theorem
 \ref{mo} the cut off function $G(\eta)$ was
 removed, because the eliminated part is a
 pseudodifferential operator which is
   bounded on any $L^p$, when
 $\tilde{m}\leq0$.

\par As an alternative, one could also use
the failure of the boundedness on $M^p$
combined with Lemma \ref{lloc} $(ii)$,
but the above approach seems a little
bit shorter. \par We observe that the
operators $\tilde{A}$ and
$\tilde{A}^\ast$ above allowed us to
prove in \cite{cordero-nicola-rodino}
the sharpness of the threshold
$-d|1/2-1/p|$ for FIOs acting on local
$\Fur L^p$ spaces.

\section*{Acknowledgements}
The authors would like to thank Sandro
Coriasco and Michael Ruzhansky for
fruitful conversations and comments.


\begin{thebibliography}{10}

\bibitem{beurling53} A. Beurling and H. Helson.
 \newblock Fourier transforms with bounded powers.
 \newblock {\it Math. Scand.},
1:120--126, 1953.


\bibitem{benyi} A. B\'enyi, K. Gr\"ochenig,
 K.A. Okoudjou and L.G. Rogers.
 \newblock Unimodular Fourier multipliers for
 modulation spaces.
 \newblock {\it J. Funct. Anal.},
246(2):366-384, 2007.


\bibitem{cappiello} M.
Cappiello.
\newblock Fourier integral
operators of infinite order and
applications to SG-Hyperbolic
equations.
\newblock{\it Tsukuba J. Math.}, 28:
311--361, 2004.


\bibitem{concetti-toft}
F.~Concetti and J.~Toft.
\newblock
{Schatten-von Neumann properties for
Fourier integral operators with
non-smooth symbols, I}.
\newblock
{\it Ark. Mat.}, to appear.

\bibitem{cordero2}
E.~Cordero and F. Nicola.
\newblock {Metaplectic representation
on
 Wiener amalgam spaces and applications
to the Schr\"odinger equation}.
\newblock {\em J. Funct.
Anal.}, 254:506--534, 2008.
\bibitem{cordero-sharp} E.
Cordero, F. Nicola.
\newblock Sharpness of some properties of Wiener amalgam and modulation
spaces.
\newblock {\em Preprint}, March 2008. Available at
ArXiv:0803.3140v1.
\bibitem{fio1}
E.~Cordero, F. Nicola and L. Rodino.
\newblock {T}ime-frequency {A}nalysis of {F}ourier {I}ntegral {O}perators.
\newblock {\em Preprint}, October 2007. Available at
ArXiv:0710.3652v1.
\bibitem{cordero-nicola-rodino}
E.~Cordero, F. Nicola and L. Rodino.
\newblock Boundedness of Fourier Integral
Operators on $\cF L^p$ spaces. {\em
Trans. Amer. Math. Soc.}, to appear.
Available at ArXiv:0801.1444.
\bibitem{co2} H. O.Cordes. {\it The thechnique of pseudodifferential operators}, Cambridge University Press, 1995.
\bibitem{cordes}
H. O. Cordes.
\newblock {\em Precisely predictable {D}irac observables}.
\newblock Fundamental Theories of Physics,
154, Springer Dordrecht, 2007.

\bibitem{coriasco} S.
Coriasco.
\newblock Fourier integral
operators in SG classes I. Composition
theorems and action on SG Sobolev
spaces.
\newblock{\it Rend. Sem. Mat.
Univ. Pol. Torino}, 57: 249--302, 1999.

\bibitem{coriasco1} S.
Coriasco.
\newblock Fourier Integral
Operators in SG classes II. Application
to SG Hyperbolic Cauchy Problems.
\newblock{\it Ann. Univ. Ferrara-Sez. VII-Sc. Mat.}, XLIV:
81--122, 1998.
\bibitem{coriasco2} S.
Coriasco and P. Panarese.
\newblock Fourier Integral
Operators Defined by Classical Symbols
with Exit Behaviour.
\newblock{\it Math. Nachr.}, 242:61--78, 2002.

\bibitem{coriascoruz} S.
Coriasco and M. Ruzhansky.
\newblock Global $ L^p$ estimates for Fourier Integral Operators.
\newblock {\it
In preparation}.
\bibitem{fefferman} C.\ Fefferman.
\newblock The uncertainty
principle. \newblock{\it Bull. Amer.
Math. Soc.,} 9:129--205, 1983.
\bibitem{F1} H.~G.~Feichtinger.
\newblock Modulation spaces on locally
compact abelian groups.
\newblock {\em Technical Report, University Vienna, 1983,}
 and also in
\newblock {\em Wavelets and Their Applications},
M. Krishna, R. Radha, S. Thangavelu,
editors,
\newblock Allied Publishers, 99--140, 2003.
\bibitem{fe89-1}
H.~G. Feichtinger.
\newblock {A}tomic characterizations of modulation spaces through {G}abor-type
representations.
\newblock In {\em {P}roc. {C}onf. {C}onstructive {F}unction {T}heory}, {\em {R}ocky {M}ountain {J}. {M}ath.}, 19:113--126, 1989.
\bibitem{feichtinger90}
H.~G. Feichtinger.
\newblock Generalized amalgams, with applications to {F}ourier transform.
\newblock {\em Canad. J. Math.}, 42(3):395--409, 1990.
\bibitem{fg97jfa}
H.~G. Feichtinger and K.~Gr{\"o}chenig.
\newblock Gabor frames and time-frequency analysis of distributions.
\newblock {\em J. Funct. Anal.}, 146(2):464--495, 1997.
\bibitem{Feichtinger-grochenig89}
H. G. Feichtinger and K. Gr\"ochenig.
\newblock Banach spaces related to integrable group representations and their atomic decompositions, II.
\newblock{\it
Monatsh. Math.}, 108:129--148, 1989.
\bibitem{folland89}
G.~B. Folland.
\newblock {\em Harmonic analysis in phase space}.
\newblock Princeton Univ. Press, Princeton, NJ, 1989.
\bibitem{book} K. Gr\"ochenig. {\it Foundation of Time-Frequency Analysis}. Birkh\"auser, Boston MA, 2001.
\bibitem{GL01}
K.~Gr{\"o}chenig and M.~Leinert.
\newblock Wiener's lemma for twisted convolution and {G}abor frames.
\newblock {\em J. Amer. Math. Soc.}, 17:1--18, 2004.

\bibitem{hormander0} L. H\"ormander.
\newblock Fourier integral
operators I.
\newblock {\em Acta Math.},
\newblock 127:79--183, 1971.

\bibitem{hormander}
 L.~H\"{o}rmander.
\newblock {\it The Analysis of Linear Partial
Differential Operators, Vol. III, IV}.
Springer-Verlag, 1985.


\bibitem{lebedev94}
V.~Lebedev and A.~Olevski\v {\i}.
\newblock $C^1$ changes of variable: Beurling-Helson type theorem and H\"ormander conjecture
on Fourier multipliers.
\newblock {\em Geom. Funct. Anal.}, 4(2):213--235, 1994.
\bibitem{mascarello-rodino} M. Mascarello and L. Rodino.
\newblock Partial Differential Equations with Multiple Characteristics.
Berlin, Akad. Verl., 1997

\bibitem{kasso07} K.A. Okoudjou.
\newblock A Beurling-Helson type theorem for
 modulation spaces. {\it J. Funct. Spaces Appl.},
 to appear.
 Available at http://www.math.umd.edu/$\sim$kasso/publications.html.

\bibitem{krantz}
S. G. Krantz and H. R. Parks.
\newblock {\it The implicit function theorem}.
\newblock Birkh\"auser Boston Inc, Boston, 2002.

\bibitem{mel1} {R. B. Melrose}. {\rm Spectral and Scattering theory of the Laplacian on asymptotically Euclidean spaces}. In {\it Spectral and Scattering theory}, M. Ikawa, ed., Marcel Dekker, 85--130, 1994.
\bibitem{mel2} {R. B. Melrose}. {\it Geometric scattering theory}. Cambridge University Press, 1995.

\bibitem{parenti72}
C. Parenti.
\newblock Operatori pseudo-differenziali in {$\R\sp{n}$} e applicazioni.
   \newblock {\em Ann. Mat. Pura Appl. (4)},
   93:359--389, 1972.
\bibitem{phongstein97}
D. H. Phong and E. M. Stein.
\newblock The {N}ewton polyhedron and oscillatory integral operators.
   \newblock {\em Acta Math.},
   179(1):105--152, 1997.

\bibitem{roch-tachiz98}
R. Rochberg and K. Tachizawa.
\newblock Pseudodifferential operators, {G}abor frames, and local
           trigonometric bases.
   \newblock {\em Gabor analysis and algorithms, Appl. Numer. Harmon. Anal.},
   171--192, Birkh\"auser Boston, 1998.

\bibitem{ruz1}
M. Ruzhansky.
\newblock On the
sharpness of Seeger-Sogge-Stein orders.
\newblock {\em Hokkaido Math. J.}, 28:357--362,
1999.

\bibitem{ruz2}
M. V. Ruzhansky.
\newblock
Singularities of affine fibrations in
the regularity theory of Fourier
integral operators. \newblock {\em
Russian Math. Surveys}, 55:93--161,
2000.

\bibitem{ruzhsugimoto}
  M. Ruzhansky and M. Sugimoto.
  \newblock Global $L^2$-boundedness theorems for a class of Fourier integral operators.
  \newblock{\em Comm. Partial Differential Equations},
   31(4-6):547-569, 2006.
\bibitem{sc1} E. Schrohe. {Spaces of weighted symbols
 and weighted Sobolev spaces on manifolds}.
 {\em LNM} 1256:360--377, Springer-Verlag, Berlin-Heidelberg, 1987.

\bibitem{sch1} B.-W. Schulze, {\it Boundary Value Problems and Singular Pseudo-differential Operators.} J.Wiley \& Sons, Chichester, New York, $1998$.

\bibitem{seegersoggestein}
A. Seeger, C. D. Sogge and E. M. Stein.
\newblock {Regularity properties of {F}ourier integral operators}.
\newblock {\em Ann. of Math. $(2)$}, 134(2):231--251, 1991.

\bibitem{sogge93}
C. D. Sogge.
\newblock {\it Fourier Integral in Classical Analysis}.
\newblock Cambridge Tracts in Math. \#105, Cambridge Univ. Press, 1993.

\bibitem{stein93}
E. M. Stein.
\newblock {\it Harmonic analysis}.
\newblock Princeton University Press, Princeton, 1993.


\bibitem{sugimototomita}
M. Sugimoto and N. Tomita.
\newblock {{T}he dilation property of modulation spaces and their inclusion relation with Besov spaces}.
\newblock {\em J. Funct. Anal.}, 248(1):79--106, 2007.

\bibitem{Toftweight}
J.~Toft.
\newblock Continuity properties for modulation spaces, with applications
              to pseudo-differential calculus. {II}.
\newblock {\em Ann. Global Anal. Geom.}, 26(1):73--106, 2004.

\bibitem{treves}
F. Tr{e}ves
\newblock
{\it Introduction to pseudodifferential
operators and Fourier integral
 operators, Vol. I, II.}
 \newblock
 Plenum Publ. Corp., New York, 1980.


 \end{thebibliography}
\end{document}